\documentclass[final,letterpaper]{siamart250211}
\usepackage{stmaryrd,amssymb}
\usepackage{amsmath}
\usepackage{tabularx,multirow}
\usepackage{picture,slashbox}
\usepackage{graphicx,caption}
\usepackage{epsfig,epsf,psfrag,epstopdf,psfig}
\usepackage[hang,center]{subfigure}
\usepackage{color,bm, mathtools}
\usepackage{comment}
\usepackage{algorithm,algorithmic}
\usepackage{url}
\usepackage{booktabs}
\usepackage{relsize}

\newcommand{\gm}{\gamma}

\newcommand\bx{\mathbf{x}}
\newcommand\bh{\mathbf{h}}
\newcommand\by{\mathbf{y}}
\newcommand\bz{\mathbf{z}}

\newcommand\bR{\mathbb{R}}

\newcommand{\Og}{\Omega}
\newcommand{\fl}[2]{\frac{#1}{#2}}

\newcommand{\lela}{\left \langle}
  \newcommand{\rira}{\right \rangle}

\renewcommand{\theequation}{\arabic{section}.\arabic{equation}}
\newcommand{\be}{\begin{equation}}
\newcommand{\ee}{\end{equation}}
\newcommand{\ba}{\begin{array}}
\newcommand{\ea}{\end{array}}
\newcommand{\bal}{\begin{aligned}}
\newcommand{\eal}{\end{aligned}}
\def\bea{\begin{eqnarray}}
\def\eea{\end{eqnarray}}
\def \beas{\begin{eqnarray*}}
\def \eeas{\end{eqnarray*}}
\newtheorem{exmp}{Example}[section]
\newtheorem{remark}{Remark}

\DeclareMathOperator{\sign}{sign}

\newcommand\threeitems[3]{%
\item#1%
\hspace{15pt}%
\hspace{\labelsep}#2
\hspace{10pt}
\hspace{\labelsep}#3
}

\newcommand\twoitems[2]{%
\item#1%
\hspace{15pt}%
\hspace{\labelsep}#2
}

\SetSymbolFont{stmry}{bold}{U}{stmry}{m}{n}

\renewcommand{\thefigure}{\arabic{figure}}
\renewcommand{\thetable}{\arabic{table}}

\title{An efficient and accurate numerical method for computing the ground states of three-dimensional rotating dipolar Bose-Einstein condensates under strongly anisotropic trap
\thanks{This work was supported by the National Key R\&D Program of  China (No. 2024YFA1012803) (Q. Tang, S. Zhang and Y. Zhang),
the Institutional Research Fund from Sichuan University (No. 2020SCUNL110) (Q. Tang and S. Zhang), the basic research fund of Tianjin University  (No. 2025XJ21-0010) (Y. Zhang),
the National Natural Science Foundation of China No.12271400 (S. Zhang and Y. Zhang) and No. 11871418, 11971120 (H. Wang), and Yunnan Fundamental Research Projects (No. 202101AS070044) (H. Wang).}}
\author{Qinglin Tang\thanks{College of Mathematics, Sichuan University, No.24 South Section 1, Yihuan Road,
Chengdu, China, 610065. ({\tt qinglin\_tang@scu.edu.cn}, {\tt shaobo\_zhang@scu.edu.cn}).}\and
Hanquan Wang\thanks{School of Statistics and Mathematics, Yunnan University of Finance and Economics,
 Kunming, Yunnan, China, 650221.
 ({\tt wang\_hanquan@hotmail.com}).} \and
Shaobo Zhang \footnotemark[2]  \and
Yong Zhang\thanks{Corresponding author. Center for Applied Mathematics and KL-AAGDM, Tianjin University, Tianjin, China, 300072.
({\tt Zhang\_Yong@tju.edu.cn}).}
}

\begin{document}
\maketitle

\begin{abstract}
In this article, we propose an efficient and spectrally accurate numerical method to compute the ground states
of three-dimensional (3D) rotating dipolar Bose-Einstein condensates (BEC) under strongly anisotropic trapping potentials.
The kernel singularity, convolution non-locality and density anisotropy together complicate the dipolar potential evaluation.
The fast rotation mechanism not only  induces a complicated energy landscape with many local minima,
but also creates a large number of vortices in the condensates.
Such factors collectively make the ground state computation challenging in terms of convergence, accuracy and efficiency, especially for 3D anisotropic systems.
Coupled with Fourier spectral discretization, we proposed a preconditioned conjugate gradient method (PCG) by integrating the anisotropic truncated kernel method (ATKM)
for the dipolar potential evaluation.
An adaptive step size control strategy is designed and ATKM allows for a spectral accuracy without introducing any extra anisotropy-dependent memory requirement or computational time.
Our algorithm is spectrally accurate, highly efficient and memory-economic.
Extensive numerical results are presented to confirm the accuracy and efficiency, together with
applications to study impacts of the model parameters on critical rotational frequency, energies and chemical potential.
Furthermore, these simulations reveal additional novel ground state patterns, such as bent vortices.

\end{abstract}


\pagestyle{myheadings}\thispagestyle{plain}
\markboth{Q. Tang, H. Wang, S. Zhang and Y. Zhang}
{three-dimensional anisotropic rotating dipolar Bose-Einstein condensates}

\section{Introduction}

Since the first successful realization of Bose-Einstein condensates (BEC) in 1995 \cite{Anderson, Bradley, Davis},
the field has witnessed numerous significant advancements \cite{Bar2008, Gaunt13, woo05}.
Notably, the introduction of the rotational frame gives rise to vortices \cite{AD2003, Yi2018BdG}, while the long-range dipolar potential leads to intriguing phenomena such as magnetostriction \cite{Magnetostriction2007} and the formation of quantum droplets \cite{droplet2016} that can self-organize into supersolid phases \cite{Casotti2024, Poli2025, supersolid}.
Recently, microwave shielding technique has been utilized to suppress the two- and three-body losses in sodium-cesium molecules,
and, for the first time, a dipolar molecular BEC has been achieved, which offers a clean and highly tunable system for strongly dipolar quantum matter and quantum simulation \cite{dipmol2024}.
In addition, by breaking the symmetry of the confining potential, anisotropic traps serve as a versatile platform for investigating a wide range of fundamental mechanisms and properties,
from the morphology of the quantum gas and anisotropic expansions \cite{Gaunt13} to the vortex arrangement \cite{Ok2002} and quasiparticle excitations \cite{woo05}.
A remarkable finding is that \textsl{bent vortices} \cite{AD2003, GP64,GP63} were observed in the elongated rotating BEC with strong local interactions, which further inspired the exploration of vortex lines.
Given these advances, the study of rotating dipolar BEC under strongly anisotropic traps  attracts great interest in the physics community,
promising the emergence of further novel quantum phenomena.


At temperature $T$ much smaller than the critical temperature $T_c$,  the property of three-dimensional rotating dipolar BEC can be described by the complex-valued wave function $\psi(\bx,t)$
whose evolution is governed by the dimensionless Gross-Pitaevskii equation (GPE) \cite{BC2013} as follows
\bea\label{DipGPE0}
&& i\partial_t \psi({\bx}, t) = \left[-\fl{1}{2}\nabla^2 + V({\bf x}) + \beta |\psi|^2 +
\lambda\, \Phi(\bx,t) -\Omega  L_z\right]\psi(\bx,t), \\
\label{dipole-poten0}
&&\Phi(\bx,t)=\left(U_{\rm dip}\ast |\psi|^2\right)(\bx,t),\\ 
\label{ini-con0}
&&\psi(\bx,0)=\psi_0(\bx), 
\eea
where $t$ is the time variable and $\ast$ represents the convolution operator
with respect to the spatial variable $\bx=(x,y,z)^T\in  {\mathbb R}^3$.
The rotation term $L_z = -i (x \partial_y - y\partial_x) = -i \partial_{\theta}$ is the
$z$-component of the angular momentum and $\Omega$ represents the rotational frequency.
The dimensionless constant $\beta$ describes the strength of short-range two-body
interactions in a condensate,
and $V(\bx)$ is a given real-valued external trapping potential
which is determined by the system under investigation.  In most BEC experiments,
a harmonic potential is chosen to trap the condensate, i.e.,\
\be\label{Vpoten}
V(\bx) = \fl{1}{2}(\gm_x^2x^2 + \gm_y^2y^2 + \gm_z^2z^2),
\ee
where $\gm_x>0$, $\gm_y>0$ and $\gm_z>0$ are  dimensionless constants proportional to the
trapping frequencies in the $x$-, $y$- and $z$-direction, respectively.
The trapping potential is isotropic if $\gm_x=\gm_y=\gm_z$, and is anisotropic otherwise.
The system is considered strongly anisotropic when maximum of the  trapping frequency ratio is great than 1,  i.e., $\max_{i,j=x,y,z}{\{\gm_i/\gm_j\}}\gg1$.
Moreover, $\lambda$ is a constant characterizing the strength of dipolar potential and
$U_{\rm dip}(\bx)$ is the dipole-dipole interaction kernel defined as
\be\label{DDI-3D}
U_{\rm dip}(\bx)= \fl{3}{4\pi |\bx|^3}\left[ 1-\fl{3 (\bx\cdot\bm{n})^2}{|\bx|^2} \right]
=-\delta(\bx)-3\,\partial_{\bm{n}\bm{n}}   \left( \fl{1}{4\pi|\bx|} \right),    \quad \ \ \bx\in{\mathbb R}^3,
\ee with dipole orientation $\bm{n} = (n_1, n_2, n_3)^T$ being a given unit vector, i.e., $|\bm{n}|=\sqrt{n_1^2+n_2^2+n_3^2}=1$,  $\partial_{\bm{n}}=\bm{n}\cdot \nabla$ and
$\partial_{\bm{n}\bm{n}}=\partial_{\bm{n}}(\partial_{\bm{n}})$.
In fact, for smooth densities $\rho(\bx) := |\psi(\bx)|^2$, the dipolar potential can be reformulated via 3D Coulomb potential
as follows
\bea\label{DDI2Cou3D}
\Phi(\bx) =-\rho - 3 \;\; \frac{1}{4\pi|\bx|}  \ast (\partial_{\bm{n}\bm{n}}\rho)
:=-\rho - 3 \; [U  \ast (\partial_{\bm{n}\bm{n}}\rho)],\quad  U(\bx) = \frac{1}{4\pi |\bx|},
\eea
where $U(\bx)$ is the three-dimensional Coulomb interaction kernel.
The time-dependent GPE (\ref{DipGPE0})--(\ref{ini-con0}) conserves two important quantities:  the
{\it total mass} (or normalization) of the wave function
\be
\label{total_mass}
N (\psi(\cdot,t)):=\|\psi(\cdot, t)\|_{L^2}^2=\int_{{\mathbb R}^3} |\psi(\bx,t)|^2 {\rm d}\bx \equiv1,\;\; t\geq 0,
\ee
and the {\it energy per particle}
\bea\label{energy}
~\qquad E(\psi(\cdot, t))\!:=\!\int_{{\mathbb R}^3}\!\left(\!\fl{1}{2}|\nabla\psi|^2 \!+\! V({\bf x})|\psi|^2\!+\! \frac{\beta}{2}  |\psi|^4 \!+ \! \fl{\lambda}{2}\Phi |\psi|^2
	\!-\!\Omega \bar{\psi} L_z\psi \!\right)\! {\rm d}\bx
\!\equiv\! {E} (\psi(\cdot, 0)),
\eea
where $\bar{\psi}$ is the conjugate of $\psi$.

To find a stationary solution of Eq. \eqref{DipGPE0}, we write $\psi(\bx,t)=\phi(\bx)e^{-i\mu t}$, where $\mu$ is the chemical potential of the condensate and $\phi(\bx)$ is the stationary state. Plugging it into \eqref{DipGPE0}, we obtain
\bea\label{EL}
\mu~\phi(\bx)=\left[-\fl{1}{2}\nabla^2 + V({\bf x}) + \beta |\phi|^2 +
\lambda\,\left(U_{\rm dip}\ast |\phi|^2\right) -\Omega  L_z\right]\phi(\bx).
\eea
This is indeed a nonlinear eigenvalue problem under a normalization constraint $N(\phi)=1$, and any eigenvalue $\mu$ can
be computed from its corresponding eigenfunction $\phi$ by
\beas
\mu&=&\int_{{\mathbb R}^3}\left(\fl{1}{2}|\nabla\phi|^2 + V({\bf x})|\phi|^2+ \beta|\phi|^4 + \lambda\left(U_{\rm dip}\ast |\phi|^2\right) |\phi|^2
	-\Omega \bar{\phi} L_z\phi \right) {\rm d}\bx \\
   &=&E(\phi)+\int_{{\mathbb R}^3}\left(\frac{\beta}{2} |\phi|^4 + \fl{\lambda}{2}\left(U_{\rm dip}\ast |\phi|^2\right) |\phi|^2\right) {\rm d}\bx.
\eeas
Equivalently, the stationary states $\phi$ are critical points of the energy functional $E(\phi)$ over the unit sphere $\mathcal{S}:=\{\phi(\bx) ~\! | ~\! \|\phi\|^2_{L^2}\!=\!1, E(\phi)\!<\!\infty\}$.
The ground state $\phi_g$ is defined as the minimizer of the energy functional $E(\phi)$:
\be\label{groundDef}
\phi_g =\arg \min_{\phi\in \mathcal{S}} \; E(\phi).
\ee

The existence of ground state in rotating dipolar BEC has been well-studied theoretically \cite{BC2013}.
When the trapping potential is strongly anisotropic, e.g., pancake-shaped or cigar-shaped,
the three-dimensional system can be approximated by a two/one-dimensional model \cite{BJNY}.
However, a gap persists between such low-dimensional approximations and physical systems,
especially in regimes of intermediate anisotropy.
Importantly, some physical phenomena, such as vortex lines in rotating BEC, can only be observed in three-dimensional space.
It is therefore essential to resort to the original three-dimensional model, despite the inevitably high computational costs involved.
To date, few numerical methods have been developed to compute the ground state of rotating dipolar BEC in strongly anisotropic regimes.

For numerical study on the ground state of anisotropic rotating dipolar BEC, there are three challenges:
(i) The strongly anisotropic trap introduces vastly different scales of variation along each direction,
demanding sufficiently high resolution in all directions to capture fine structures.
(ii) The singular kernel and the anisotropic density pose significant challenges to evaluating the dipolar potential with both accuracy and efficiency.
(iii) Fast rotation induces the formation of quantum vortices and creates an intricate energy landscape with many local minima \cite{shu2024preconditioned},
which complicates the computation of the ground state.

\begin{figure}[h!]
\subfigure{
\psfig{figure=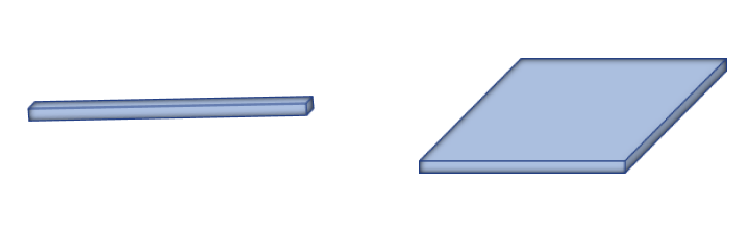,height=2.0cm,width=12.5cm,angle=0}
}
\caption{Schematic diagrams of cigar- and pancake-shape domain $\mathcal{D}_{L\bm{\xi}}$ with $\bm{\xi}=(\frac{1}{10},\frac{1}{10},1)^T$ (left) and $\bm{\xi}=(1,1,\frac{1}{10})^T$ (right).
}\label{fig_domain}
\end{figure}

In presence of a confining anisotropic potential $V(\bx)$, the wave function $\psi$ (density $\rho$) is smooth, fast-decaying and anisotropic. Thus,
it suffices to truncate the whole space to an anisotropic rectangular domain
\bea \label{domain}
\mathcal{D}_{L\bm{\xi}}:=\prod\nolimits_{\alpha} \, [-L \xi_\alpha, L \xi_\alpha],  \qquad ~0< \xi_\alpha \leq 1, \quad \quad \alpha=x, y, z,
\eea
with $\xi_\alpha\!=\!1$ on the longest side of $\mathcal{D}_{L\bm{\xi}}$.
The truncated domain is then discretized with uniform mesh grid in each spatial direction.
The anisotropy vector and strength are defined as
\bea\label{Def-AnisoStr}
\bm{\xi} = (\xi_x, \xi_y, \xi_z)^T, \qquad \xi_f=1/(\xi_x \xi_y \xi_z)
\eea
respectively. For readers' convenience, we present graphical illustrations of the anisotropic computation domains
with $\bm{\xi}=(\frac{1}{10},\frac{1}{10},1)^T$ (cigar-shape)
and $\bm{\xi}=(1,1,\frac{1}{10})^T$ (pancake-shape) in Fig. \ref{fig_domain}.
The wave function (density) is well approximated by the Fourier spectral method \cite{BC2013,SpectralBkShen}, which helps achieve spectral accuracy with $\mathcal{O}(N\!\log\! N)$ efficiency thanks to the discrete Fast Fourier Transform (FFT).
To capture fine structures such as vortex lattices,
the mesh size should be small enough, which usually leads to a large system
and consequently severe memory and computational costs.

Due to kernel singularity, convolution non-locality and density anisotropy,
it is challenging to compute the dipolar potential with spectral accuracy while maintaining the FFT-like efficiency.
During the last decade, several fast spectral methods have been proposed,
including the NonUniform Fast Fourier Transform method (NUFFT) \cite{Nufft}, Gaussian-Summation method (GauSum) \cite{GauSum}, Kernel Truncation Method (KTM) \cite{KTM} and  Anisotropic Truncated Kernel Method (ATKM)\cite{atkm}.
All these methods achieve similar spectral accuracy but differ in terms of efficiency, particularly in cases where the density is strongly anisotropic.
For NUFFT and KTM, the memory requirement and efficiency performance
depend linearly on the anisotropy strength, therefore, they can hardly handle the strongly anisotropic case.
The GauSum method splits the potential into one regular and singular integral, and costs more computation time than ATKM.
In this article, we adopt ATKM, which is introduced by Greengard {\sl et al.}~\cite{atkm}
to especially deal with the anisotropic density via rectangular instead of circular extension of the computation domain $\mathcal{D}_{L\bm{\xi}}$ \eqref{domain}.
The memory requirement and efficiency of ATKM are independent of the anisotropy strength, which guarantees excellent performance even when the density function is strongly anisotropic.


For the computation of the ground state, there have been many methods proposed based on imaginary-time evolution \cite{BC2013, BD_SISC, BCEsine, SP_NUFFT}.
For fast rotating BEC, these methods often exhibit slow convergence, making it impractical to obtain the ground state within an affordable computational time.
Alternatively, the ground state can also be obtained by solving a constrained nonlinear eigenvalue problem,
which actually corresponds to an optimization problem defined on a Riemannian manifold \cite{Altmann2021, gseigen}.
It is quite natural to apply the state-of-art Riemannian optimization techniques \cite{Absil2009, Boumal2023} therein. For example,
a preconditioned conjugate gradient algorithm (PCG), proposed \cite{PcgNLSE17} for the rotating BEC, has greatly accelerated the convergence and significantly improved the efficiency.
The ground state of rotating dipolar BEC has been investigated by combing PCG with KTM \cite{PcgCiCP2018} for normal isotropic trapping potentials.
For the BEC under strongly anisotropic traps, the aforementioned method is hardly practical due to the huge memory burden
and long computational time brought in by the usage of KTM for the dipolar potential evaluation.

\

Regarding all the above aspects, in this article, we aim to propose a preconditioned conjugate gradient algorithm by combining ATKM for dipolar potential evaluation.
By designing an appropriate step size control strategy and addressing the challenge arising from integration of PCG and ATKM,
our method achieves spectral accuracy with high efficiency, and more importantly, the memory requirement is independent of the anisotropy strength.
Then, we study bent vortices in the  dipolar BEC under different parameter setups (Subsection \ref{bent}), and these results confirm the importance and necessity of our algorithm for handling such strongly anisotropic dipolar systems.
In addition, using our algorithm, we identify novel ground state patterns and comprehensively investigate the impacts of model parameters on the critical rotational frequency and energies.

The rest of the paper is organized as follows:
In Section \ref{sec:Alg}, we introduce the spatial discretization by Fourier spectral method and a brief review of the ATKM for dipolar potential evaluation,
and then propose an accurate and efficient numerical algorithm by integrating PCG with ATKM.
In Section \ref{ResultsSec}, we present extensive numerical results to confirm the accuracy and efficiency,
along with applications to investigate the impacts of model parameters on  critical rotational frequency, energies and chemical potential.
Moreover, we identify some ground state patterns and study the bent vortices in dipolar BEC.
Finally, some conclusions are drawn in Section \ref{sec:Conclu}.

\section{Numerical algorithm}\label{sec:Alg}
Here we first utilize the Fourier spectral method \cite{SpectralBkShen} to discretize the  wave function $\phi$ on domain $\mathcal{D}_{L\bm{\xi}}$ \eqref{domain}.
Then, we present a detailed introduction of ATKM for evaluating dipolar potentials, and propose an accurate and efficient numerical method to compute the ground state  of  rotating dipolar BEC
under anisotropic traps.

\subsection{Spatial discretization}\label{spatial_four}
The computation domain $\mathcal{D}_{L\bm{\xi}}$ is discretized along each spatial direction with mesh size
$h_{\alpha}=2L\xi_{\alpha}/N_{{\alpha}}$, where $N_{\alpha} (\alpha=x,y,z)$ are even positive integers,
and  the total grid number is denoted as $N=N_xN_yN_z$.
For simplicity, we denote the mesh size vector $\bh:=(h_x,h_y,h_z)$, grid number vector $\bm{N}:=(N_x, N_y, N_z)$
and mesh grids $\bx_{\bm{j}}:=(x_{j_x},y_{j_y},z_{j_z})$.
The set of  mesh grids is defined as
$$\mathcal{T}:=\{\bx_{\bm{j}}=(-L\xi_x+{j_x}h_x,-L\xi_y+{j_y}h_y,-L\xi_z+{j_z}h_z), ~ \bm{j}\in \mathcal{I}_{\bm{N}}\},$$
with index set
$$\mathcal{I}_{\bm{N}}:=\{\bm{j}=(j_x,j_y,j_z)\in \mathbb{N}^3~ \big|~ 0 \leq j_\alpha \leq N_\alpha-1, ~\alpha = x,y,z \}.$$
We denote $f_{\bm{j}}$ as the numerical approximation of $f(\bx_{\bm{j}}),~\bx_{\bm{j}} \in \mathcal{T}$ and $f$ as
the vector with components $\{f_{\bm{j}}\}_{\bm{j}\in \mathcal{I}_{\bm{N}}}$.
Moreover, we introduce the Fourier index set
$$\Lambda_{\bm{N}}:=\{\bm{k}:=(k_x,k_y,k_z)\in \mathbb{Z}^3~ | -N_\alpha/2 \leq k_\alpha\leq N_\alpha/2-1, ~\alpha = x,y,z\},$$
and define the Fourier basis function as
\beas\label{Fourbasis}
W_{\bm{k}}^L(\bx):= e^{i \nu_{k_{x}}(x+L\xi_{x})}~e^{i \nu_{k_{y}}(y+L\xi_{y})}~e^{i \nu_{k_{z}}(z+L\xi_{z})},
\quad \nu_{k_{\alpha}}=\pi k_{\alpha}/(L\xi_{\alpha}),\quad \bm{k}\in \Lambda_{\bm{N}}.
\eeas
The Fourier spectral approximation of wave function $\phi$ on domain $\mathcal{D}_{L\bm{\xi}}$ reads as
\bea\label{Fourphi}
\phi(\bx)\approx
\sum_{\bm{k}\in \Lambda_{\bm{N}}}
\widehat{\phi}_{\bm{k}}~W_{\bm{k}}^L(\bx), \qquad \bx \in \mathcal{D}_{L\bm{\xi}},
\eea
where $\widehat{\phi}_{\bm{k}}$, the discrete Fourier coefficients, are given by
\bea\label{FourphiCo}
\widehat{\phi}_{\bm{k}}:=
\frac{1}{N}\sum_{\bm{j}\in \mathcal{I}_{\bm{N}}} \phi_{\bm{j}}~\overline{W_{\bm{k}}^L}(\bx_{\bm{j}}).
\eea
The above summation can be accelerated by FFT within
$\mathcal{O}(N\log N)$ operations.
In fact, the approximation \eqref{Fourphi} is actually an interpolation function at mesh grid points $\bx_{\bm{j}}\in\mathcal{T}$.
We  approximate  $\partial_x\phi$, $\Delta \phi$ and $L_z\phi$ on mesh $\mathcal{T}$ by differentiating the corresponding finite Fourier series \eqref{Fourphi} as follows:
\bea\label{Laplacian}
(\partial_{x}\phi)_{\bm{j}}&\approx&
(\llbracket{\partial_{x}}\rrbracket\phi)_{\bm{j}}=\sum_{\bm{k}\in \Lambda_{\bm{N}}}
(i\nu_{k_x})~\widehat{\phi}_{\bm{k}}~W_{\bm{k}}^L(\bx_{\bm{j}}), \notag\\[-0.2em]
(\Delta\phi)_{\bm{j}}&\approx& (\llbracket \Delta \rrbracket\phi)_{\bm{j}}~
=\sum_{\bm{k}\in \Lambda_{\bm{N}}}
\left(-\sum_{\alpha=x,y,z}(\nu_{k_{\alpha}})^2\right)\widehat{\phi}_{\bm{k}}~W_{\bm{k}}^L(\bx_{\bm{j}}),\\[0.4em]
(L_z\phi)_{\bm{j}}&\approx& (\llbracket L_z \rrbracket\phi)_{\bm{j}}
=-i(x\llbracket{\partial_{y}}\rrbracket\phi-y\llbracket{\partial_{x}}\rrbracket\phi)_{\bm{j}}.\notag
\eea

\subsection{Anisotropic truncated kernel method}
\label{algSect}
To evaluate the dipolar potential, it is clear that direct discretization of the convolution integral might diverge when the source and target points coincide,
so special treatment is required to deal with such singularity that is caused by the singular but integrable convolution kernel.
Numerical schemes based on direct integral discretization, whose computation complexity is about $\mathcal{O}(N^2)$ flops, are often too time-consuming and become even impractical
for high-dimensional problems since the total grid number $N$ will easily soar to millions or billions.
Moreover, under a strongly anisotropic trap, the density function exhibits a strong spatial anisotropy. That is to say,
the span of its compact support in one or two spatial directions is much shorter than those in the other directions.
Such anisotropic density will aggravate the computation difficulty.
Here we choose the fast spectral solver ATKM \cite{atkm}, which is particularly designed to evaluate convolution-type nonlocal potentials with anisotropic densities.
As shown in Eq. \eqref{DDI2Cou3D}, the dipolar potential can be rewritten as a summation of the density and a Coulomb potential generated by $\partial_{\bm{n}\bm{n}}\rho(\bx)$.

Since the density $\rho$ is smooth, anisotropic and decays exponentially fast, it is reasonable to assume
that $\rho$ is compactly supported in an anisotropic rectangular domain $\mathcal{D}_{L\bm{\xi}}$.
The potential on $\mathcal{D}_{L\bm{\xi}}$ now reads
\be\label{iden_b2}
\bal
\Phi(\bx) &=  \int_{\bR^3} U(\by)  \rho(\bx-\by) {\rm d} \by
= \int_{\bx + \mathcal{D}_{L\bm{\xi}}} U(\by)  \rho(\bx-\by) {\rm d} \by, \\
&=\int_{\mathcal{D}_{2L\bm{\xi}}} U(\by)  \rho(\bx-\by) {\rm d} \by,
\quad \;\;\bx \in \mathcal{D}_{L\bm{\xi}}.
\eal
\ee
The third equation holds because, for any $\bx \in \mathcal{D}_{L\bm{\xi}},\, \by \in \mathcal{D}_{2 L\bm{\xi}}\backslash(\bx+\mathcal{D}_{L\bm{\xi}})$,
we have $\bx-\by \notin \mathcal{D}_{L\bm{\xi}}$ and thus $\rho(\bx-\by)=0$.
To integrate Eq. \eqref{iden_b2}, we needs to approximate the density $\rho$ on domain $\mathcal{D}_{3L\bm{\xi}}$,
because $\bx-\by \in \mathcal{D}_{3L\bm{\xi}}$ for any $\bx\in\mathcal{D}_{L\bm{\xi}}$ and $\by \in \mathcal{D}_{2L\bm{\xi}}$.
It is natural to extend the density to $\mathcal{D}_{3L\bm{\xi}}$ by zero-padding and apply the Fourier spectral method therein.
Fortunately, it is sufficient to approximate the density on a \textit{twofold} domain $\mathcal{D}_{2L\bm{\xi}}$ thanks to the periodicity of Fourier series,
and we refer to \cite{atkm,GauSum} for more details.

The density $\rho$ is well resolved by the following Fourier spectral method
\be\label{FourSeriB3}
\rho(\bz) \approx \sum_{\bm{k}\in\Lambda_{2\bm{N}}} \widehat{\rho}_{\bm{k}}
~W_{\bm{k}}^{2L}(\bz)
, \qquad \bz \in \mathcal{D}_{2L\bm{\xi}},
\ee
where the discrete Fourier coefficients $\widehat{\rho}_{\bm{k}}$ are given below
\be\label{fouriercoeff}
\widehat{\rho}_{\bm{k}}:=
\frac{1}{8N}\sum_{\bm{j}\in \mathcal{I}_{2\bm{N}}} \rho_{\bm{j}}~\overline{W_{\bm{k}}^{2L}}(\bz_{\bm{j}}).
\ee
Plugging \eqref{FourSeriB3} into \eqref{iden_b2}, we obtain
\be \label{U_TrunFour}
\bal
\Phi(\bx) 
\approx \sum_{\bm{k}\in\Lambda_{2\bm{N}}} \widehat{U}_{\bm{k}} \;\widehat{\rho}_{\bm{k}}\;
~W_{\bm{k}}^{2L}(\bx), \quad \;\;\bx \in \mathcal{D}_{L\bm{\xi}},
\eal
\ee
where $\widehat{U}_{\bm{k}}$, the Fourier coefficients of kernel $U(\by)$ on domain $\mathcal{D}_{2L\bm{\xi}}$, are defined below
\be\label{FUk}
\widehat{U}_{\bm{k}}:=\int_{\mathcal{D}_{2L\bm{\xi}}} U(\by) ~\overline{W_{\bm{k}}^{2L}}(\by)~{\rm d }\by.
\ee
The above singular Fourier integral is actually well-defined, and Greengard \textit{et al.} proposed an efficient and accurate algorithm
to evaluate it based on a black-box sum-of-Gaussians approximation of the symmetric kernel $U(\by)$ \cite{atkm}.
By taking advantage of the separable structure of Gaussian functions,
the Fourier coefficients $\widehat{U}_{\bm{k}}$ can be calculated within $\mathcal{O}(N)$ operations.
Clearly, once $\widehat{U}_{\bm{k}}$ is available, the nonlocal potential computation \eqref{U_TrunFour} can be accomplished
in three simple steps: (1) compute the discrete Fourier coefficient $\widehat\rho_{\bm{k}}$ \eqref{fouriercoeff} via FFT;
(2)compute $\widehat{U}_{\bm{k}}\,\widehat \rho_{\bm{k}}$ by pointwise function multiplication and (3) compute $\Phi(\bx_j)$ via inverse FFT.
The total computational costs amount to $\mathcal{O}(2^3N\log(2^3N))$.

As mentioned above, the density anisotropy property is taken into account at the very beginning by taking an anisotropic rectangular domain
and carrying out a twofold anisotropic extension.
Once the computation of $\widehat{U}_{\bm{k}}$ \eqref{FUk}, which is a Fourier integral defined on the anisotropic domain too, is available, the potential is computed via a pair of FFT/iFFT on vectors of twofold length in each spatial direction.
Thus, both the computational time and memory requirement are \textsl{independent} of the anisotropy strength $\xi_f$ \eqref{Def-AnisoStr}.
While, for KTM, the total optimal zero-padding factor linearly depends on the anisotropy strength \cite{optKTM}, and
the precomputation of Fourier coefficients suffers badly from prohibitive memory requirement and computational costs in high dimensions.

\begin{remark}
To provide a vivid illustration of the comparison between KTM and ATKM,
we consider a double-precision computation with an anisotropy vector $\bm{\xi}=(\frac{1}{10},\,\frac{1}{10},\,1)^T$.
Here we compare the memory requirement of Fourier coefficients $\widehat{U}_{\bm{k}}$.
Given a uniform mesh with grid number $N\!\!\!=\!\!\!256^3$,
it requires only $256^3\times2^3\times8/2^{30}$Gb~$=$~\textbf{1} Gb memory for ATKM.
In contrast, the memory costs of KTM, whose total optimal zero-padding factor \cite{KTM} amounts to
$432$, 
soars as much as \textbf{54} Gb, which is far beyond most personal computers' capacity.
In other words, the reduction factor on memory requirement is $1- \frac{8}{432}\approx \bm{98\%}$.
\end{remark}

\begin{remark}
In simulation, the nonlocal potentials are often called multiple times with the same numerical setups, that is, both the domain and mesh grid are fixed.
Therefore, the evaluation of $\widehat{U}_{\bm{k}}$ is computed \textsl{once for all} in the precomputation step.
\end{remark}

\subsection{Ground state computation}
Preconditioned conjugate gradient methods have been established for computing the ground states of rotating BEC with/without dipolar potential \cite{PcgNLSE17, PcgCiCP2018, zhang2024}.
While their efficiency is well-established, existing approaches are designed under either isotropic or weakly anisotropic traps.
For the most challenging strongly anisotropic regime,
we develop a preconditioned conjugate gradient method that incorporates ATKM for dipolar potential evaluation to compute the ground state,
addressing the difficulties existing in the subtle vortex-structure ground state caused by fast rotation, the normalization constraint and dipolar potential under a strongly anisotropic trap.

\

Based on the spatial discretization by Fourier pseudo-spectral method in subsection \ref{spatial_four}, we discretize
the energy $E(\phi)$ on mesh grid $\mathcal{T}$ as follows
\bea
\qquad~~ E(\phi)\approx \mathcal{E}(\phi)=\left\langle -\frac{1}{2}\llbracket{\Delta}\rrbracket\phi+V\odot\phi+\frac{\beta}{2}|\phi|^2\odot\phi
+\frac{\lambda}{2}\llbracket{\Phi\rrbracket}\odot\phi-\Og\llbracket{L_z}\rrbracket\phi, \;\phi \right\rangle,
\eea
where $\llbracket{\Phi\rrbracket}$ is the numerical approximation of $\Phi$, $\odot$ denotes Hadamard product and
\vspace{-0.2cm}
\bea\label{innerPro}
\langle U, V\rangle := \Re\Biggl( h_x h_y h_z\sum_{\bm{j}\in \mathcal{I}_{\bm{N}}} U_{\bm{j}}\overline{V}_{\bm{j}} \Biggl),
\qquad \forall \, U, V\in \mathbb{C}^{N},
\eea
is the inner product defined on $\mathbb{C}^{N}$ with $\Re(\cdot)$ representing the real part of a complex number.
The minimization problem \eqref{groundDef} is then discretized into the finite dimensional minimization problem
\bea
\label{mp}
\phi_g=\arg\min_{\phi\in \mathbb{S}}\mathcal{E}(\phi), \quad \mathbb{S}:=\{\phi \in \mathbb{C}^{N}|~ \|\phi\|=1 \},
\eea
with norm $\|U\|=\langle U, U\rangle^{\frac{1}{2}}$.

The (Euclidean) gradient of energy $\mathcal{E}(\phi)$ on $\mathbb{C}^{N}$ is
\bea
&&\nabla \mathcal{E}(\phi)=2H_{\phi}\phi, \\
&&H_{\phi}\phi= -\frac{1}{2}\llbracket{\Delta}\rrbracket\phi+V\odot\phi+\beta|\phi|^2\odot\phi
+\lambda\llbracket{\Phi\rrbracket}\odot\phi-\Og\llbracket{L_z}\rrbracket\phi. \notag
\eea
The tangent space at a point $\phi \in \mathbb{S}$ is defined as $T_{\phi}\mathbb{S}:=\{f \in  \mathbb{C}^{N}~|~ \langle f,\phi\rangle=0\}$.
And the orthogonal projection $J_{\phi}(f)$ from $\mathbb{C}^{N}$ onto the tangent space $T_{\phi}\mathbb{S}$ is given by
\begin{equation}
J_{\phi}(f):=f-\langle f,\phi\rangle\,\phi.
\end{equation}
The Riemannian gradient of energy on $\mathbb{S}$ is
\bea
\nabla_{\mathbb{S}} \mathcal{E}(\phi)=J_{\phi}(\nabla \mathcal{E}(\phi)) =2(H_{\phi}\phi-\mu_{\phi} \phi)
\quad \mbox{with} ~~ \mu_{\phi}=\langle H_{\phi}\phi, \phi\rangle.
\eea
The first-order necessary optimality condition states that at a minimum $\phi \in \mathbb{S}$, the gradient of the energy vanishes on the constraint manifold $\mathbb{S}$, i.e., $\nabla_{\mathbb{S}} \mathcal{E}(\phi)=0$. This condition is indeed the Euler-Lagrange equation associated with problem \eqref{mp}:
\bea
H_{\phi}\phi = \mu_{\phi}\,\phi.
\eea
Further, the preconditioned conjugate gradient method on manifold $\mathbb{S}$ is proposed.

\

To numerically solve the minimization problem \eqref{mp}, we denote the abstract grid function $f$ of the $n$-th step as $f_n$.
Moving from $\phi_n$ along a descent direction $p_n$ on the manifold $\mathbb{S}$, we use the retraction
\bea\label{retraction}
\phi_{n+1} = R_{\phi_n}(t_np_n):= \cos(t_{n}\|p_n\|) \phi_{n} + \sin(t_n\|p_{n}\|) \,\,p_{n}/\|p_{n}\|,
\eea
with $p_n \in T_{\phi_n}\mathbb{S}$ and step size $t_n > 0$ (to be determined later).
Regarding the direction $p_n$, we first introduce the tangent bundle of $\mathbb{S}$ as $T\mathbb{S}:=\cup_{\phi\in \mathbb{S}} T_{\phi}\mathbb{S}$.
The transport operator $\mathcal{V}: T\mathbb{S} \times T\mathbb{S}\longrightarrow T\mathbb{S}$, given explicitly below
\beas
\mathcal{V}_u(v):=J_{R_{\phi}(u)}(v)=v-\langle v, R_{\phi}(u)\rangle R_{\phi}(u), \quad \forall\, u\in T_{\phi}\mathbb{S}, ~v\in T_{\varphi}\mathbb{S},
\eeas
maps vector $v\in T_\varphi\mathbb{S}$ at point $\varphi$ to tangent space at point $R_\phi(u)\in \mathbb{S}$.
Specifically, for $p_{n\!-\!1}\in T_{\phi_{n\!-\!1}}\mathbb{S}$ and $v\in T_{\phi_{n\!-\!1}}\mathbb{S}$, we have $\mathcal{V}_{t_{n\!-\!1}p_{n\!-\!1}}(v)=J_{\phi_n}(v)$.
Then the direction of conjugate gradient method is chosen as
\bea
~\qquad~~ p_n=
\begin{cases}
\!J_{\phi_n}(-Pr_n), \qquad\qquad\qquad \qquad\quad \qquad\qquad\qquad \qquad\quad~\qquad\quad  n=0,\\
\!J_{\phi_n}(-Pr_n)\!+\!\beta^{cg}_n \mathcal{V}_{t_{n\!-\!1}p_{n\!-\!1}}(p_{n\!-\!1})
\!=\!J_{\phi_n}(-Pr_n)\!+\!\beta^{cg}_nJ_{\phi_n}( p_{n\!-\!1}),
 ~n>0,\\
\end{cases}
\eea
with the residual $r_n=\nabla_{\mathbb{S}} \mathcal{E}(\phi_n)/2=H_{\phi_n}\phi_n-\mu_{\phi_n}\phi_n$.
The parameter $\beta^{cg}_n$ is given by Polak-Ribi\`ere choice \cite{Absil2009}
\bea
&&\beta^{cg}_n=\max(\beta^{PR}_n, 0), \\
&& \beta^{PR}_n=\langle r_n - J_{\phi_n}(r_{n-1}),
J_{\phi_n}(Pr_{n})\rangle/\langle r_{n-1}, J_{\phi_{n-1}}(Pr_{n-1})\rangle. \notag
\eea
The preconditioner $P$ is a symmetric positive definite operator to help accelerate the convergence of iteration,
which is given by \cite{PcgCiCP2018}
\bea
P=P_V^{1/2}P_{\Delta}P_V^{1/2},
\eea
where $P_\Delta$ and $P_V$ are respectively read as
\bea
&&P_\Delta=(\alpha_{\Delta}-\llbracket{\Delta\rrbracket}/2)^{-1}, \\
&&P_V=\left( \alpha_V+ V +\beta|\phi_n|^2+|\lambda|(1+\sign(\llbracket{\Phi_n\rrbracket}))\llbracket{\Phi_n\rrbracket}/2 \right)^{-1},\notag
\eea
with $\alpha_\Delta$ and $\alpha_V$ being positive shift constants, and here we take
\bea
~\qquad \alpha_\Delta\!=\!\alpha_V\!=\!\langle-\frac{1}{2}\llbracket{\Delta\rrbracket}\phi_n
+ V\!\odot\!\phi_n
+\beta|\phi_n|^2\!\odot\!\phi_n, \phi_n\rangle
\!+\!|\lambda\langle\llbracket{\Phi_n\rrbracket}\!\odot\!\phi_n, \phi_n\rangle| >0.
\eea
This combined preconditioner achieves nearly optimal convergence performance, whose iteration number is independent of the computation domain and mesh size. For more details, we refer readers to \cite{PcgNLSE17, PcgCiCP2018}.

To find $t_n$ in Eq. \eqref{retraction}, we solve the nonlinear minimization problem
\bea\label{Opt_Theta}
t_{n} =\arg \min_{t>0}\;Q(t), \quad\quad Q(t):=\mathcal{E}\left(R_{\phi_n}(t\,p_n)\right),
\eea
by utilizing a standard one-dimensional minimization algorithm.
Since the exact line search is time-consuming, we introduce the following approximation method to obtain a simple and much cheaper approximation of $t_n$.
Now that the retraction \eqref{retraction} is of second order accuracy, we expand  $Q(t)$ as follows
\bea\label{Energy_Expend}
\qquad ~~~~
Q(t)\approx \mathcal{E}(\phi_n)+ 2t \lela H_{\phi_n}\phi_n, \;p_n  \rira
 +t^2 \big( \lela H_{\phi_n} p_n,\;p_n \rira
  + \lela g_n,\,p_n\rira-\mu_{\phi_n}\|p_n\|^2\big),
\eea
with $g_n= 2\phi_n \left(\beta \Re(\bar{\phi}_n p_n)+\lambda\,U\ast \Re(\bar{\phi}_n p_n) \right)$.
Note that $\lela H_{\phi_n}\phi_n, \;p_n  \rira$ may not be negative,
so the direction $p_n$ might not be a descent direction and the energy might not decrease.
To select an appropriate step size $t_n$, we follow this procedure:
(i) We check whether the coefficient of the first-order term is negative or not.
If so, we choose the conjugate gradient direction and keep $\beta^{cg}_n$ unchanged;
otherwise, we employ the preconditioned steepest descent (PSD) method \cite{PcgNLSE17}, i.e., we set  $\beta^{cg}_n=0$ so to ensure that $\lela H_{\phi_n}\phi_n, \;p_n  \rira<0$
given that the preconditioner $P$ is positive definite.
(ii) We check further whether the coefficient of the second-order term is positive or not.
If so, minimizing  \eqref{Energy_Expend} with respect to $t$ yields
\be
\label{opt_step}
t_n^{\rm app}=-\lela H_{\phi_n}\phi_n, \;p_n  \rira /
\left( \lela H_{\phi_n} p_n,\;p_n \rira + \lela g_n,\,p_n\rira-\mu_{\phi_n}\|p_n\|^2 \right).
\ee
If not, we set the step size to a default value, i.e., $t_n = 0.3$.
(iii) We then check whether the energy of $\phi_{n+1}$ decreases or not.
If so, we accept the step size; otherwise, we decrease it by half until the energy decreases.
The above process guarantees the energy decay property.

For the convenience of readers,
the proposed method for solving the minimization problem \eqref{groundDef} is summarized in {\bf Algorithm}  \ref{alg:CG}, which we refer to as PCG-ATKM.
In this paper, we adopt the following stopping criterion
\begin{equation}\label{stop_maxnorm}
\|\bm{\phi}_{n+1}-\bm{\phi}_n\|_{l^\infty} \leq \varepsilon.
\end{equation}

\begin{algorithm}[H]                  
\caption{\quad PCG-ATKM }          
\label{alg:CG}                           
\begin{algorithmic}                    
  \STATE Precompute the coefficients $\widehat{U}_{\bm{k}}$.
  \vspace{0.1cm}
   \STATE  Set the iteration counter $n=0$ and the initial data $\phi_0$.
  \vspace{0.25cm}
    \WHILE{\;\it not converged\;}    \vspace{0.05cm}
    \STATE {\bf Step 1:} Compute $\Phi_n$  via ATKM.
  \vspace{0.1cm}
   \STATE {\bf Step 2:} $r_{n} = H_{\phi_{n}}\phi_{n} - \mu_{\phi_n} \phi_{n}$ with $\mu_{\phi_n} = \lela H_{\phi_{n}}\phi_{n}, \phi_{n}\rira$.
   \vspace{0.17cm}
   \STATE {\bf Step 3:} $p_n=
\begin{cases}
J_{\phi_n}(-Pr_n), \qquad\qquad\qquad\qquad   n=0,\\
J_{\phi_n}(-Pr_n)+\beta^{cg}_nJ_{\phi_n}( p_{n-1}), \quad \mbox{others}.\\
\end{cases}$
 \vspace{0.17cm}
   \STATE {\bf Step 4:} Correct the direction $p_n$ and choose step size $t_n$:
   \vspace{0.05cm}
   \begin{itemize}
\item []~~~~~{\scriptsize$\bullet$} If $\lela H_{\phi_n}\phi_n, \;p_n  \rira<0$, keep $p_n$ unchanged;\\
    ~~~~~~~otherwise, set $p_n=J_{\phi_n}(-Pr_n)$.    \vspace{0.05cm}
\item []~~~~~{\scriptsize$\bullet$} If $\lela H_{\phi_n} p_n,\;p_n \rira
  + \lela g_n,\,p_n\rira-\mu_{\phi_n}\|p_n\|^2>0$, take $t_n$ as \eqref{opt_step};\\
    ~~~~~~~otherwise, set $t_n=0.3$.    \vspace{0.05cm}
\item []~~~~~{\scriptsize$\bullet$}
If $\mathcal{E}(\phi_{n+1})\geq\mathcal{E}(\phi_{n})$, decrease $t_n$ by half until the energy decreases.
\end{itemize}
\vspace{0.1cm}
   \STATE {\bf Step 5:} $\phi_{n+1} = \cos(t_{n}\|p_n\|) \phi_{n} + \sin(t_n\|p_{n}\|){p_{n}}/{\|p_{n}\|}$.
   \vspace{0.1cm}
   \STATE {\bf Step 6:} $n = n+1$.    \vspace{0.05cm}
    \ENDWHILE
\end{algorithmic}
\end{algorithm}

During each iteration, both the descent direction determination
and step size approximation involve dipolar potentials evaluations.
Given the spatial anisotropy exhibited by the wave function (density) in an anisotropic trap,
the adoption of ATKM suits best the strong anisotropic setup
in terms of memory requirement and computational cost.
Using the same uniform mesh grid on the same domain, the combination of preconditioned conjugate gradient method and ATKM is quite natural and requires no additional efforts.

Direct solving the minimization problem \eqref{groundDef} on a fixed fine mesh requires computations of high-dimensional functionals,
therefore, it is quite time-consuming even though most of the computations can be accelerated by FFT,
especially for 3D anisotropic system.
The initial guesses are often inaccurate approximations of the ground state,
which significantly postpones the convergence process, especially for cases that involve fast rotation and strong anisotropy.
To alleviate such inefficiency, we adopt a cascadic multigrid technique \cite{PcgCiCP2018}.
Numerical solution interpolated using coarser-mesh solution usually provides a better approximation than the commonly-used initial guesses,
hence a faster convergence is expected and observed indeed in computation practice.
Such cascadic multigrid approach is summarized in {\bf Algorithm} \ref{alg:CG_Multi-Grid}. Hereafter, we refer to it as PCG-ATKM-MG for convenience.

\begin{algorithm}[H]
\caption{\quad PCG-ATKM-MG}
\label{alg:CG_Multi-Grid}
\begin{algorithmic}
   \STATE   Given initial mesh grid $\mathcal{T}^i$ with grid numbers
    $N_{\alpha}^i=k_{\alpha}\times2^i$, $\alpha=x,y,z$.
    \vspace{0.07cm}
     \STATE  Given an initial data  $\phi^i_0$ on  $\mathcal{T}^i$.
   \vspace{0.25cm}
    \WHILE{\;\it not reach the finest mesh\;}
   \vspace{0.07cm}
   \STATE {\bf Step 1:} Compute the ground state $\phi^i_g$ by utilizing {\bf  Algorithm} {\rm \ref{alg:CG}}.
   \vspace{0.07cm}
   \STATE {\bf Step 2:} Set  $\phi^{i+1}_0$ by interpolating $\phi^i_g$ on the refined mesh $\mathcal{T}^{i+1}$.
   \vspace{0.07cm}
   \STATE {\bf Step 3:} $i =i+1$.
   \vspace{0.07cm}
    \ENDWHILE
\end{algorithmic}
\end{algorithm}

\begin{remark}
Under an isotropic trap, the algorithm PCG-KTM \cite{PcgCiCP2018} serves as a suitable, accurate and efficient numerical method for computing the ground states of rotating dipolar BEC.
While under an anisotropic trap, this method will require substantially large storage and more computational time,
since the nonlocal potential solver KTM employs an isotropic truncation technique.
As seen in \cite{optKTM}, the total optimal zero-padding factor of KTM increases linearly with the anisotropy strength.
This suggests that a prohibitively large memory requirement will be induced under an extremely strong anisotropic trap.
Notably, the memory requirement of PCG-ATKM-MG is \textbf{independent} of  the anisotropy strength.
\end{remark}

\begin{remark}
With appropriate modifications, the PCG-ATKM-MG can be extended to other related BEC models, such as multi-component BEC \cite{TwoCompDipRot}
and droplet dipolar BEC \cite{droplet2016}, which involves the famous Lee-Huang-Yang term \cite{lhy1957} to describe the quantum
fluctuation mechanism.
We do not  provide related details in this paper and shall leave them in some future work.
\end{remark}

\section{Numerical results} \label{ResultsSec}

In this section, we test the accuracy and efficiency of PCG-ATKM-MG,
and then  apply it to study various vortices patterns,
for example, bent vortices in rotating dipolar BEC.
In addition, we investigate the impacts of  model parameters on the critical rotational frequency, and study the variations of energies and chemical potential in details.

In the following experiments, we  consider the harmonic potential \eqref{Vpoten}, choose the rectangular domain $\mathcal{D}_{L\bm{\xi}}$ and  discretize it uniformly in each spatial direction.
We take an equal number of grid points in each spatial direction, and
the algorithm is executed with a three-level multigrid starting from $N/4^3$ total points to $N$ points.
As for the initial guesses, we have the following 10 choices \cite{PcgCiCP2018}:
\begin{description}
\threeitems{$(a)$ $\phi_a(\bx) = \frac{1}{\sqrt[4]{\pi^3}}  e^{\frac{-|\bx|^2}{2}} $,}
{$(b)$ $\phi_b(\bx) = \frac{(x + i y ) \phi_a(\bx)}{\|(x + i y ) \phi_a(\bx)\|_{L^2}} $,}
{$(\bar b)$ $\phi_{\bar{b}}(\bx)=\bar{\phi}_b(\bx)$,}
\twoitems
{$(c) ~\phi_c(\bx) = \frac{\phi_a(\bx)+\phi_b(\bx)}{\|\phi_a(\bx)+\phi_b(\bx)\|_{L^2}}$,}
{$(\bar c)~ \phi_{\bar{c}}(\bx)=\bar{\phi}_c(\bx)$,}
\twoitems
{$(d)~ \phi_d(\bx) =\frac{(1-\Omega)\phi_a(\bx)+\Omega\phi_b(\bx)}{\|(1-\Omega)\phi_a(\bx)+\Omega\phi_b(\bx)\|_{L^2}}$,}
{$(\bar d) ~\phi_{\bar{d}}(\bx)=\bar{\phi}_d(\bx)$,}
\threeitems{$(e)~\phi_e(\bx) = \frac{\Omega\phi_a(\bx)+(1-\Omega)\phi_b(\bx)}{\|\Omega\phi_a(\bx)+(1-\Omega)\phi_b(\bx)\|_{L^2}}$,}
{$(\bar{e})~\phi_{\bar{e}}(\bx)=\bar{\phi}_e(\bx)$,}
{$(f)~\phi_f(\bx) = \frac{\phi_g^{\rm TF}(\bx)}{\|\phi_g^{\rm TF}(\bx)\|_{L^2}\rule{0pt}{1.5ex}}$,}
\end{description}
where
{\small
$$
\phi_g^{\rm TF}(\bx)=\left\{
\begin{aligned}
\sqrt{(\mu_g^{\rm TF}-V(\bx))/\beta}, ~~V(\bx)<\mu_g^{\rm TF},\\
0, \qquad \qquad \qquad \qquad ~ {\mbox {\small otherwise}},
\end{aligned}
\right.
\quad \mbox{ {\small with} } ~~\mu_g^{\rm TF}=\frac{1}{2}\left(\frac{15\beta\gamma_x \gamma_y\gamma_z}{4\pi}\right)^{2/5}.$$}
In practice, we choose the stopping criterion \eqref{stop_maxnorm} with $\varepsilon=10^{-10}$ unless stated otherwise,
and select numerical solution with the lowest energy as the ground state, which is denoted as $\phi^{\bh}_g$ hereafter.
When the trap is chosen as the harmonic potential \eqref{Vpoten}, the energies associated with stationary state $\phi$ \eqref{EL} satisfy the following \textbf{virial identity} \cite{BC2013}
\be\label{vir_iden}
0=I(\phi):=2 E_{\rm kin}(\phi)  - 2 E_{\rm pot}(\phi)  + 3 ~E_{\rm int}(\phi) +3 E_{\rm dip}(\phi),
\ee
with
$$
E_{\rm kin}(\phi)=\int_{{\mathbb R}^3}\fl{1}{2}|\nabla\phi|^2 {\rm d} \bx, ~~
E_{\rm pot}(\phi)=\int_{{\mathbb R}^3}  V({\bf x})|\phi|^2 {\rm d} \bx,$$
$$
E_{\rm int}(\phi)=\int_{{\mathbb R}^3} \frac{\beta}{2}  |\phi|^4{\rm d} \bx, ~~
E_{\rm dip}(\phi)=\int_{{\mathbb R}^3} \fl{\lambda}{2}\Phi\, |\phi|^2{\rm d} \bx.
$$
The numerical approximation of the above virial identity \eqref{vir_iden} is denoted as $I^\bh$,
and we adopt its absolute value $|I^\bh|$ as the virial-identity residual so as to measure its approximation accuracy.
The algorithm is implemented in Fortran, and all reported timing results are obtained using a
Intel(R) Xeon(R) Gold6248R CPU \@ 3.00GHz with 36 MB cache with the Intel  ifort and optimization level -O3 in Ubuntu GNU/Linux. 

\subsection{Accuracy confirmation and efficiency  performance}

In this subsection, we shall present the spectral spatial accuracy and efficiency performance of our method
for both isotropic and anisotropic traps with different physical parameters.

Note that
the equality $E(\phi)=E(\kappa\phi)$ holds for any $\kappa \in \mathbb{C}$ with $|\kappa|=1$.
Consequently, the numerical errors are measured in relative $l^2$-norm as follows
\beas
{E_{2}}:=\|\phi_g^{\bh}-\kappa\,\phi_g\|_{\bh}/\|\phi_g\|_{\bh},  \qquad
\kappa=\langle\phi_g^{\bh},\phi_g\rangle_{\bh}/\langle\phi_g,\phi_g\rangle_{\bh},
\eeas
where $\langle f,g\rangle_{\bh}:=h_x h_y h_z\sum_{\bm{j}\in \mathcal{I}_{\bm{N}}} f_{\bm{j}}~\overline{g}_{\bm{j}}$ is the inner product of two complex-valued grid functions,
 the `exact' ground state $\phi_g$ is obtained numerically on a very fine mesh ($N=512^3$), and
$\phi_g^{\bh}$ is the numerical solution computed with mesh size $\bh$.

It is worthy to point out that, the energy is rotational invariant under isotropic traps, that is, $E(\phi(R_{\alpha}\bx))=E(\kappa\phi(\bx))$ holds for any rotational matrix
 $$R_{\alpha}=\begin{pmatrix}
   \cos\alpha &-\sin\alpha \\\sin\alpha & \cos\alpha
\end{pmatrix}, \qquad \forall~ \alpha \in [0,2\pi].$$
As the rotational frequency increases beyond the critical value, the ground state develops more than one vortex.
In such cases, it is inappropriate to use $E_2$ to measure the numerical approximation errors.
Thus we introduce a modified error defined as
\beas
\widetilde{E_{2}}:= \min_{\alpha} \|\phi_g^{\bh}(\!R_{\alpha}\bx)\!-\!\kappa_{\alpha}\phi_g(\bx) \|_{\bh}/\|\phi_g(\bx)\|_{\bh},
~ \kappa_{\alpha}=\langle\phi_g^{\bh}(\!R_{\alpha}\bx),\phi_g(\bx)\rangle_{\bh}/\langle\phi_g(\bx),\phi_g(\bx)\rangle_{\bh}.
\eeas
The above optimization problem is solved using the interior-point method \cite{interior}, which has been implemented in the built-in function 'fmincon' of Matlab.
The function-rotation mapping, from $\phi(\bx)$ to $\phi(R_\alpha\bx)$,
is implemented efficiently and accurately with FFT via the fast-rotation algorithm based on shear-decomposition \cite{soc-rotation-jcp,spin2Soc-WeiZhang}.

\begin{exmp} \label{TimeAccuracy} (\textbf{Accuracy and efficiency tests})
Given  dipole orientation $\bm{n} = (0,0,1)^T$, $\beta =100$ and $\lambda=0.8\beta$,
we compute the ground states for different rotational frequencies $\Omega$ in the following two cases:
\vspace{0.08cm}
\begin{itemize}
\item [] \textbf{Case I}(Isotropic trap): $\gamma_x=\gamma_y=\gamma_z=1$; \vspace{0.05cm}
\item [] \textbf{Case II}(Anisotropic trap): $\gamma_x=1$, $\gamma_y=2$, $\gamma_z=10$.
\end{itemize}\vspace{0.08cm}
\end{exmp}
We set the stopping criterion \eqref{stop_maxnorm} as $\varepsilon=10^{-12}$,
and choose domain $\mathcal{D}_{L\bm{\xi}} = [-16,16]^3$ for \textbf{Case I} and $\mathcal{D}_{L\bm{\xi}}=[-16,16]^2\times[-10,10]$ for \textbf{Case II} due to its anisotropic trap.

\begin{figure}[h!]
\subfigure{
\hspace{0.8cm}
\psfig{figure=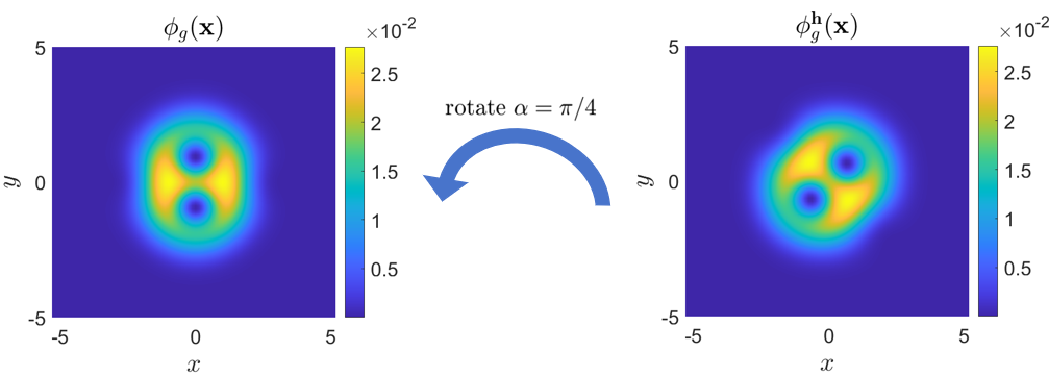,height=4.cm,width=11cm,angle=0}
}
\caption{Contour plots of $|\phi_g(x, y, z=0)|^2$ and $|\phi_g^{\bh}(x, y, z=0)|^2$ with $N=256^3$ for $\Omega=0.85$ of \textbf{Case I}  in \textbf{Example \ref{TimeAccuracy}}.
}\label{fig_rotation}
\end{figure}

\begin{figure*}[htbp]
   \centering
   \subfigure{
   \psfig{figure=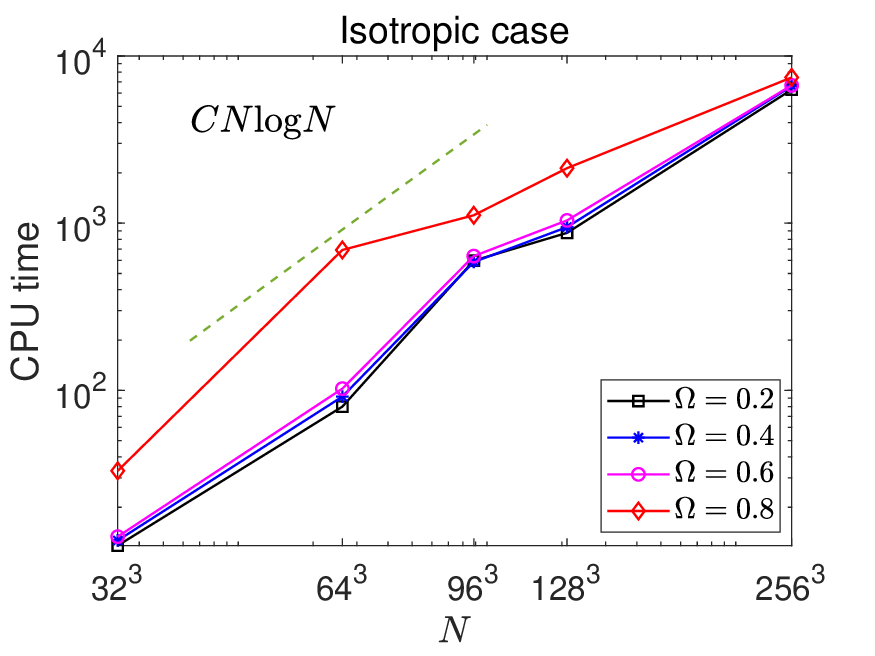,height=5.2cm,width=6.cm,angle=0}}
   \subfigure{
   \psfig{figure=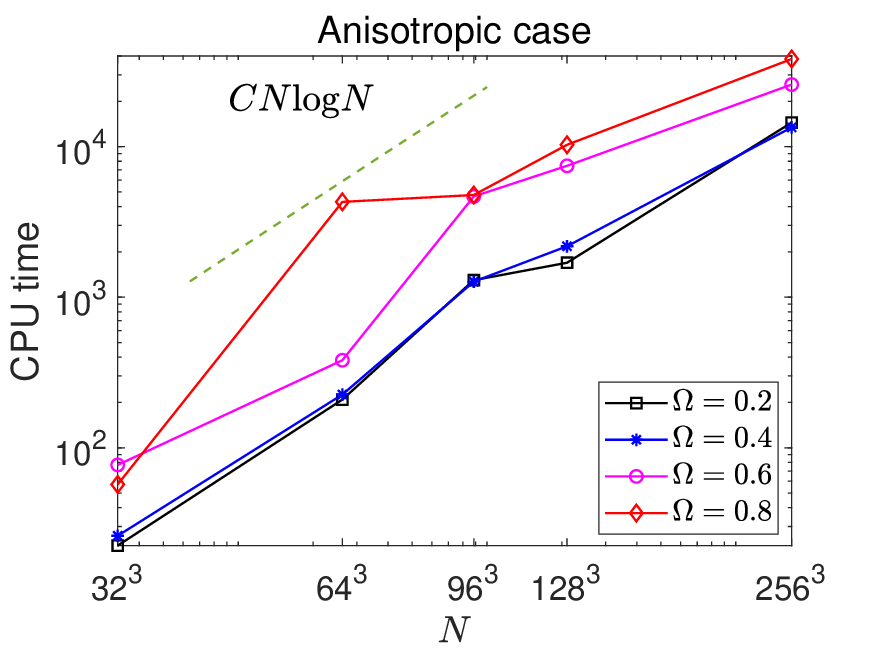,height=5.2cm,width=6.cm,angle=0}}
   \caption{Computational time for \textbf{Case I} (left) and \textbf{Case II} (right) in \textbf{Example \ref{TimeAccuracy}}. }\label{fig_efficiency}
\end{figure*}

Table \ref{tab_accuracy} shows the relative $l^2$ errors $E_2, \widetilde{E_2}$ and the virial-identity residual $|I^\bh|$ for
\textbf{Case I-II},
and Fig. \ref{fig_rotation} presents the contour plots of `exact' ground state $|\phi_g(x, y, z=0)|^2$ and numerical ground state $|\phi_g^{\bh}(x, y, z=0)|^2$ with $N=256^3$ for $\Omega=0.85$ in \textbf{Case I}. 
We can observe that errors $E_2, \widetilde{E_2}$ and  $|I^\bh|$ converge spectrally
accurate with respect to the spatial mesh size for all cases.
Fig. \ref{fig_efficiency} displays the computational time (measured in seconds), from which we know that the efficiency exhibits a typical complexity of $\mathcal{O}(N\!\log\!N)$ similar as  FFT.
This behavior agrees well with our expectation because the number of iterations in each computation
is significantly smaller than the grid number $N$.
To sum up, PCG-ATKM-MG is an accurate and efficient numerical solver for computing ground states in both isotropic and anisotropic cases.
\begin{table}[htbp]
   \centering
   \caption{Numerical errors for both isotropic and anisotropic traps in \textbf{Example \ref{TimeAccuracy}}.}\label{tab_accuracy}
   \setlength\extrarowheight{1.5pt}
   \begin{tabular}{cccccc}
      \toprule
      \multicolumn{6}{c}{\textbf{Case I:}
      $\bm{\gamma_x=\gamma_y=\gamma_z=1}$} \\
      \midrule
      $\Omega$              &  $N$           & $32^3$             &$64^3$             & $128^3$        &  $256^3$              \\
      \midrule
      \multirow{2}{*}{0.2}                   &$E_2$ & 4.1267E-03           & 2.3477E-05 & 2.1482E-09  &2.2755E-12  \\
                 &$|I^\bh|$ & 2.1743E-02      & 2.5129E-05     & 1.6822E-11  &1.3450E-12 \\
      \midrule
      \multirow{2}{*}{0.4}                  &$E_2$ &  4.2588E-03          &  2.3483E-05 &2.1482E-09 & 4.3439E-12 \\
                 &$|I^\bh|$ & 2.1810E-02      & 2.5129E-05  &1.5696E-11  &  2.0317E-13\\
      \midrule
      \multirow{2}{*}{0.6}                   &$E_2$ & 4.6460E-03           & 2.3493E-05   & 2.1482E-09  & 6.9024E-12\\
                &$|I^\bh|$ & 2.1959E-02           &  2.5129E-05  & 1.5724E-11  & 4.8439E-13 \\
      \midrule
0.8              &$E_2$ & 1.0000          &8.5676E-05   & 4.9185E-09  &  3.7105E-11\\
   (1 vortex)             &$|I^\bh|$  &  9.3167E-02     &1.2527E-04    & 2.7377E-11  &1.1491E-13\\
      \midrule
     0.85                 &$\widetilde{E_2}$ & 1.1340E-01          &7.9153E-04   & 1.4546E-07  &  3.4234E-09\\
     (2 vortices)           &$|I^\bh|$  &  1.0275E-01     &5.6180E-05    & 4.3209E-11  &6.9289E-12\\

      \midrule
      \multicolumn{6}{c}{\textbf{Case II:} $\bm{\gamma_x=1, \gamma_y=2, \gamma_z=10}$} \\
      \midrule
      $\Omega$              &  $N$           & $32^3$             &$64^3$             & $128^3$        &  $256^3$              \\
      \midrule
      \multirow{2}{*}{0.2}                   &$E_2$ & 1.0737E-01           & 2.5143E-03 & 2.4725E-06  &1.0235E-11  \\
                  &$|I^\bh|$ & 4.3461      & 5.4160E-02     & 3.7568E-06  &4.4853E-12 \\
      \midrule
      \multirow{2}{*}{0.4}                  &$E_2$ &  1.1212E-01         &  2.5217E-03 &2.4340E-06 & 3.3305E-11 \\
                 &$|I^\bh|$ & 4.4398      & 5.5158E-02  &3.5739E-06  &  2.2435E-12\\
      \midrule
      \multirow{2}{*}{0.6}                   &$E_2$ & 3.1052E-01           & 2.5985E-03   & 2.3537E-06  & 9.1651E-11\\
                &$|I^\bh|$ & 4.6578           &  5.6558E-02  & 2.9120E-06  & 3.7419E-12 \\
      \midrule
      0.8                 &$E_2$ &  7.8085E-01          & 4.7745E-02   & 8.7262E-06  &  3.5937E-10\\
      (2 vortices)       &$|I^\bh|$  &  4.5065     &2.8864E-02    & 5.6266E-06  &1.1440E-11\\
      \bottomrule
   \end{tabular}
\end{table}
%

\subsection{Bent vortices in dipolar BEC} \label{bent}
Bent vortices are known to emerge in elongated rotating BEC with strong local interactions, and their behavior in the absence of dipolar potential has been well established in previous studies \cite{AD2003, GP64,GP63}.
Building on this foundation, this subsection investigates  bent vortices in the presence of dipolar potential.
Hereafter, if not explicitly stated, we shall mix $\phi^{\bh}_g$ and $\phi_g$ just for presentation simplicity.

\begin{exmp}\label{Shape_grids} (\textbf{U-shape vortex line})
Given $\gamma_{x}=\gamma_{y}=1$, $\gamma_{z}=0.05$, $\beta=1000$, $\lambda=-300$, $\textbf{n}=(0,0,1)^T$ and $\Omega=0.75$,
we choose domain $\mathcal{D}_{L\bm{\xi}} = [-10,10]^2 \times [-100,100]$ and take grid numbers $N=128^3, 256^3, 512^3$ respectively
to compute the ground state.
\end{exmp}

\begin{figure}[h!]
\hspace{-0.3cm}
\subfigure[S-shape]{
\psfig{figure=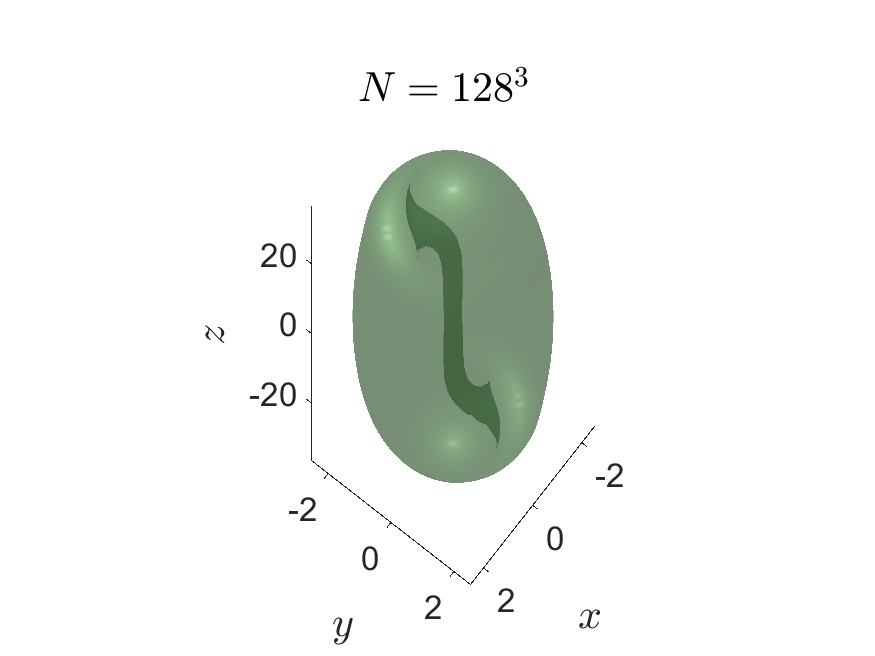,height=4.2cm,width=5.4cm,angle=0}}
\hspace{-1.825cm}
\subfigure[U-shape]{
\psfig{figure=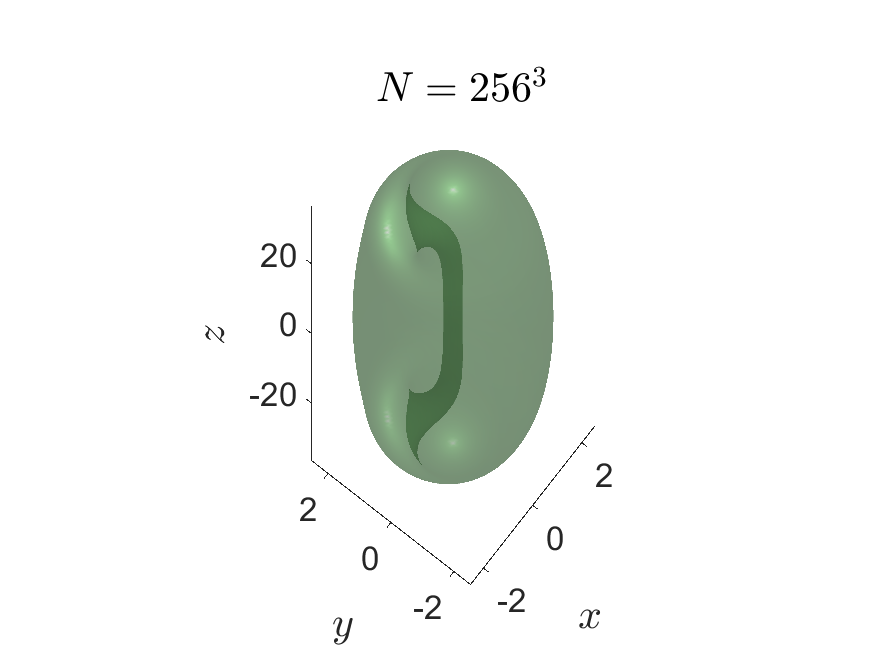,height=4.2cm,width=5.4cm,angle=0}}
\hspace{-1.895cm}
\subfigure[U-shape]{
\psfig{figure=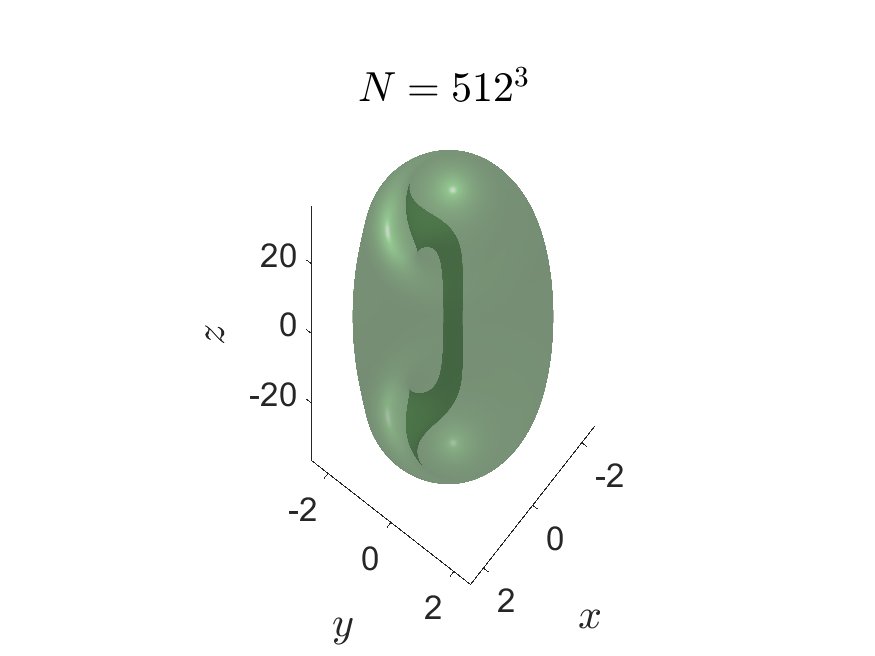,height=4.2cm,width=5.4cm,angle=0}}
\caption{Isosurface plots of $|\phi_g(\bx)|^2=3\times10^{-4}$ with different $N$ in \textbf{Example \ref{Shape_grids}}.}\label{Shape_grids_fig}
\end{figure}

Table \ref{detail_grids} presents the energies $E$, virial-identity residual $|I^\bh|$ and Hamiltonian residual $\|r_n\|$ for different $N$,
and Fig. \ref{Shape_grids_fig} shows the corresponding isosurface plots of $|\phi_g(\bx)|^2=3\times10^{-4}$.
These results indicate that the U-shape vortex state is more likely to be the ground state.
Conversely, the S-shape vortex state calculated with $N=128^3$
is not the ground state, since its associated virial-identity residual $|I^\bh|$  is relatively large with only 5-digit accuracy
( $9.0018 \times 10^{-6}$ as shown in Table \ref{detail_grids}) due to the insufficient grid number. 
The vortex structure is quite subtle and its computation is sensitive to the mesh size,
therefore, the mesh should be chosen as fine as possible, which of course poses a great numerical challenge especially for anisotropic trap
and strong dipolar potential cases.
\begin{table}[h!]
\setlength{\extrarowheight}{1.5pt}
   \centering
   \caption{Energies $E$, $|I^\bh|$ and $\|r_n\|$ for different $N$ in \textbf{Example \ref{Shape_grids}}. }\label{detail_grids}
   \begin{tabular}{m{2.7cm}  m{3.1cm} m{3.1cm} m{1.8cm}}
      \toprule
        $N$           &$E$                                      &$|I^\bh|$       & $\|r_n\| $                \\
      \midrule
         $128^3$          & 2.276210347750    & 9.0018E-06  & 5.1475E-09                \\
        $256^3$           & 2.276210390179    &2.8555E-11     &  7.8405E-09      \\
        $512^3$            & 2.276210390179   &3.0787E-12     & 7.4651E-09        \\
\bottomrule
   \end{tabular}
\end{table}

\begin{exmp}\label{DDInoDDI} (\textbf{Bent vortices under different dipolar potentials})
Given $\gamma_{x}=\gamma_{y}=1$, $\gamma_{z}=0.1$ and $\beta=9000$, we  compute the ground states  for both $\Omega=0.41$ and $\Omega=0.49$ in following cases:
\vspace{0.08cm}\begin{itemize}
\twoitems{[]\textbf{Case I}: $\lambda=0$; $\qquad\qquad\qquad\qquad ~$}{\textbf{Case II}: $\lambda=4000$, $\textbf{n}=(0,0,1)^T$;}
\vspace{0.05cm}
\twoitems{[]\textbf{Case III}: $\lambda=4000$, $\textbf{n}=(1,0,0)^T$;}
{\textbf{Case IV}: $\lambda=-4000$, $\textbf{n}=(0,0,1)^T$.}
\end{itemize}\vspace{0.08cm}
\end{exmp}

\begin{figure}[h!]
\centerline{
\psfig{figure=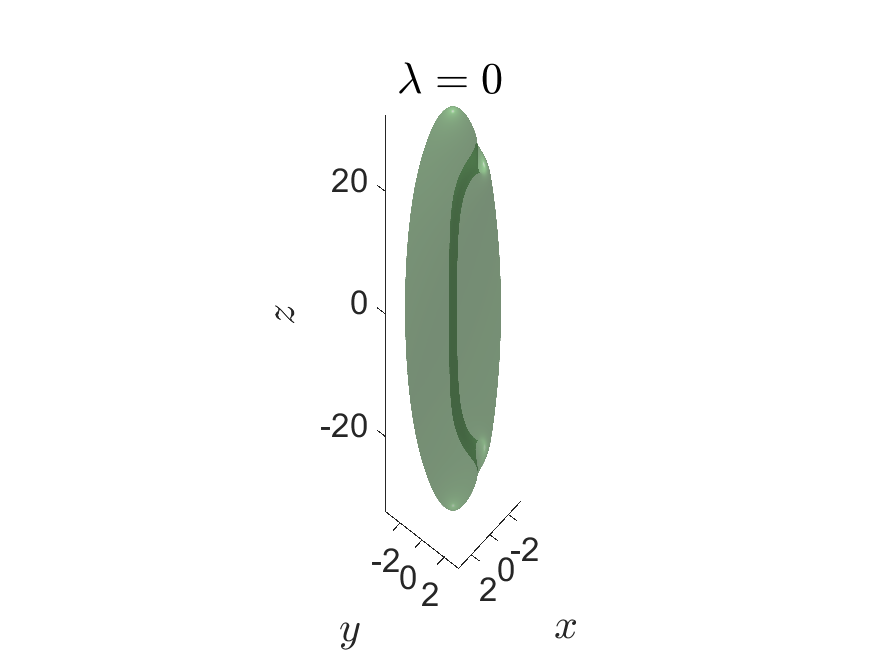,height=4.cm,width=5.5cm,angle=0}\quad
\hspace{-2.3cm}
\psfig{figure=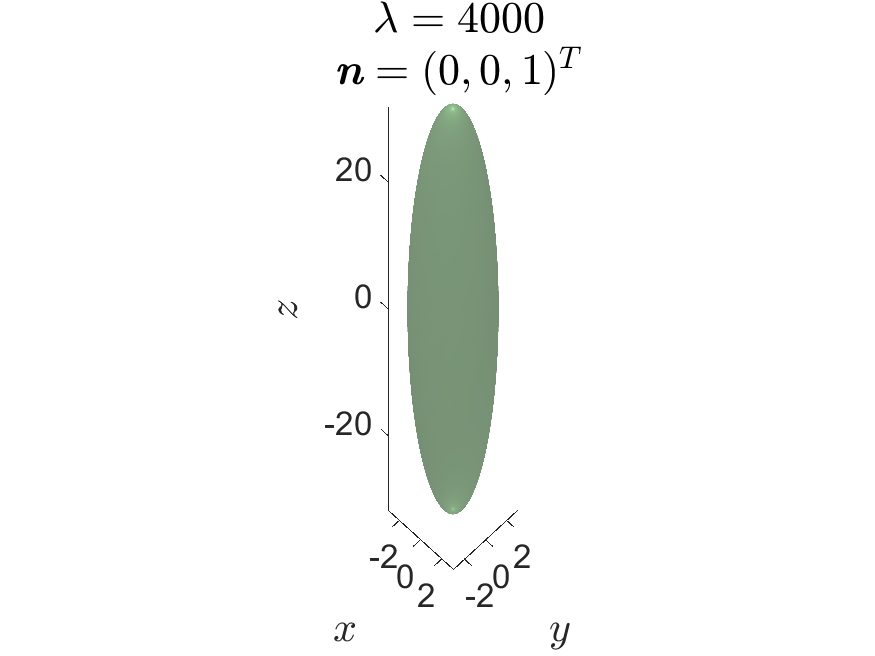,height=4.cm,width=5.5cm,angle=0}\quad
\hspace{-2.3cm}
\psfig{figure=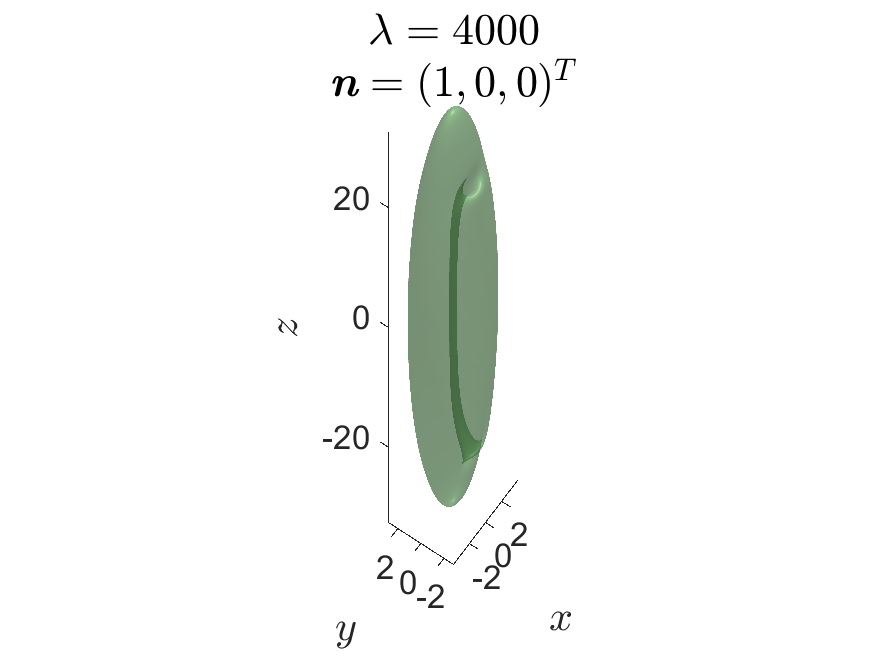,height=4.cm,width=5.5cm,angle=0}\quad
\hspace{-2.3cm}
\psfig{figure=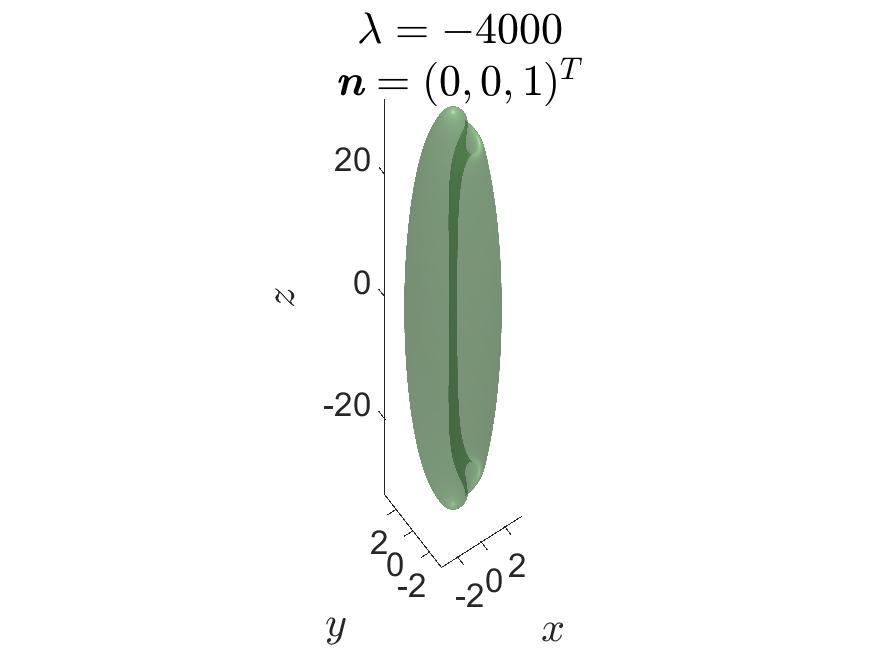,height=4.cm,width=5.5cm,angle=0}
}

\vspace{+0.6cm}

\centerline{
\psfig{figure=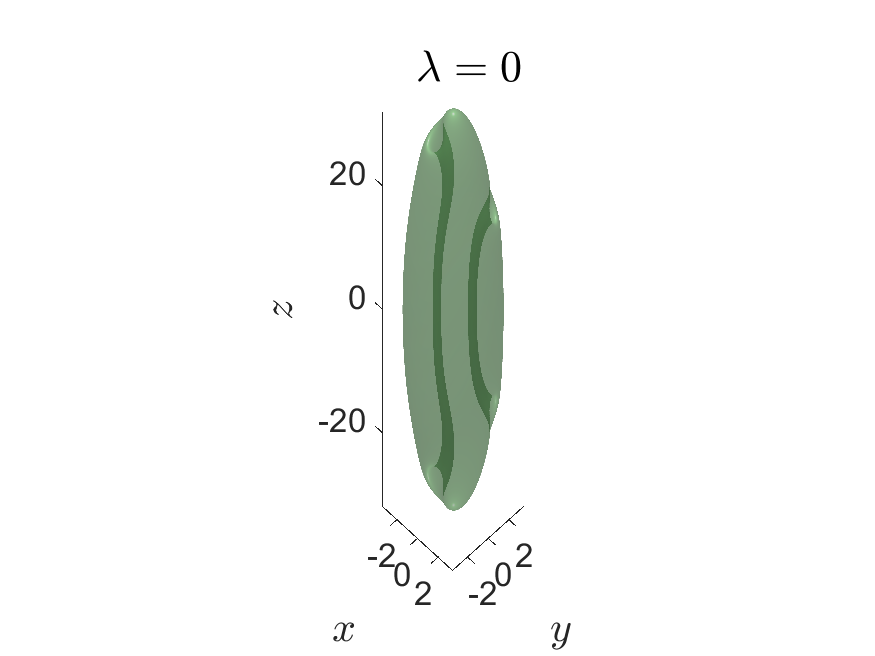,height=4.cm,width=5.5cm,angle=0}\quad
\hspace{-2.3cm}
\psfig{figure=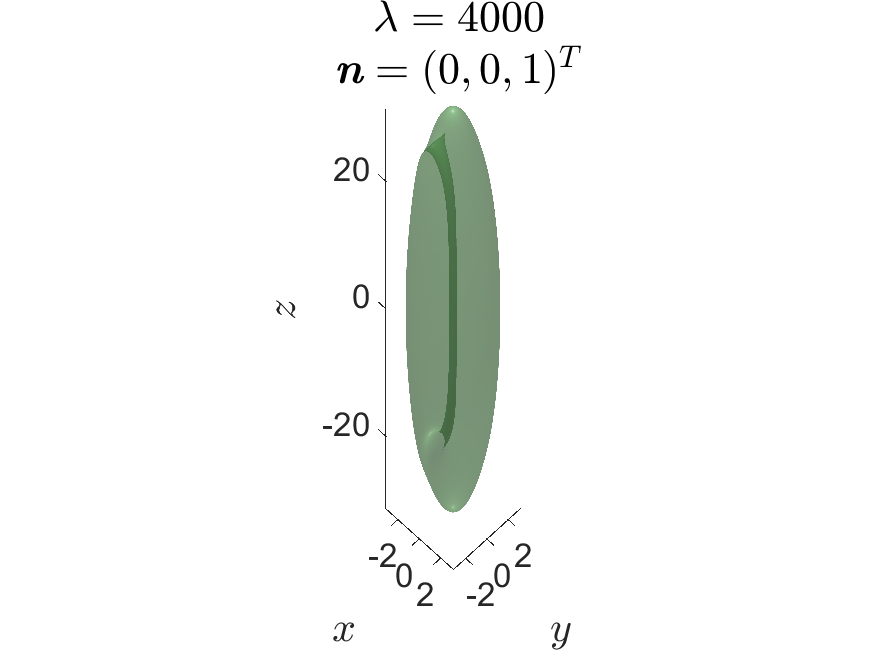,height=4.cm,width=5.5cm,angle=0}\quad
\hspace{-2.3cm}
\psfig{figure=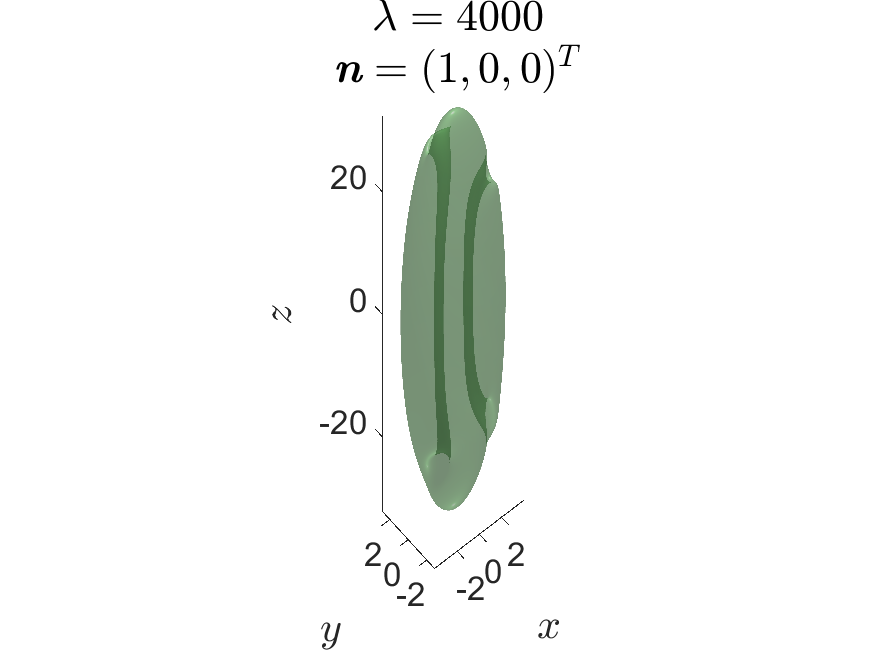,height=4.cm,width=5.5cm,angle=0}\quad
\hspace{-2.3cm}
\psfig{figure=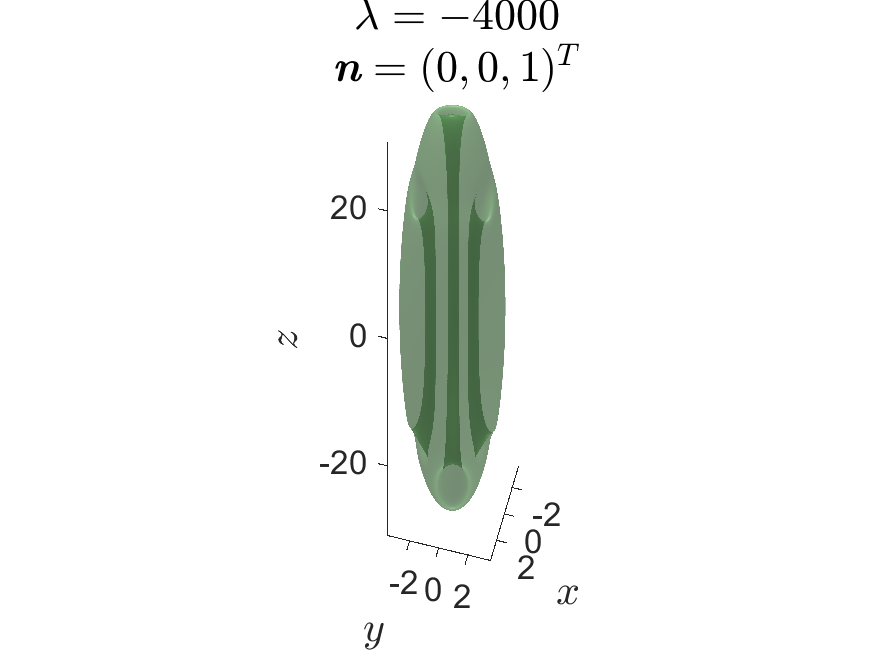,height=4.cm,width=5.5cm,angle=0}}
\caption{Isosurfaces of $|\phi_g(\bx)|^2=3\times10^{-4}$ for $\Omega=0.41$ (top row) and $\Omega=0.49$ (bottom row) in \textbf{Case I-IV} (from left to right) of \textbf{Example \ref{DDInoDDI}}.}\label{DDInoDDI_fig}
\end{figure}

The computation is done on $\mathcal{D}_{L\bm{\xi}} = [-10,10]^2 \times [-100,100]$ with $N = 384^3$.  Fig. \ref{DDInoDDI_fig} shows the isosurfaces of the densities $|\phi_g(\bx)|^2$,
from which we conclude that the dipolar potential affects vortex patterns, including number of vortices and length of the central vortex line.
In addition, for $\bm{n}=(0,0,1)^T$, there are more vortices in the ground state for negative $\lambda$ than its positive counterpart (cf. the second and fourth columns of Fig. \ref{DDInoDDI_fig}).

\subsection{Critical rotational frequency}
In rotating BEC, a vortex appears in the ground state only if the rotational frequency $\Omega$ exceeds certain
critical value $\Omega_c$. In this subsection, we shall investigate how the critical rotational frequency $\Omega_{c}$ varies with trapping frequencies of the harmonic potential \eqref{Vpoten},
local interaction strength $\beta$, dipolar potential strength $\lambda$ and the dipole orientation $\bm{n}$, respectively.
Unless stated otherwise, we take $\mathcal{D}_{L\bm{\xi}}=[-10,10]^{3}$ and $N=128^{3}$.

\begin{exmp} \label{exp_critical1}(\textbf{The impacts of external potential and local interaction})
In the absence of dipolar potential, i.e., $\lambda=0$, we investigate the variations of $\Omega_c$ in the following cases:
\vspace{0.08cm}
\begin{itemize}
\item []\textbf{Case I}: Let $\beta=2000$,  $\gamma_{x}=\gamma_{z}=1$ and vary $\gamma_{y}$;\vspace{0.05cm}
\item []\textbf{Case II}: Let $\beta=2000$, $\gamma_{x}=\gamma_{y}=1$ and vary $\gamma_{z}$;\vspace{0.05cm}
\item []\textbf{Case III}: Let $\gamma_{x}=\gamma_{y}=\gamma_{z}=1$ and vary $\beta$.
\end{itemize}\vspace{0.08cm}
\end{exmp}

\begin{figure}[h!]
\centerline{
\psfig{figure=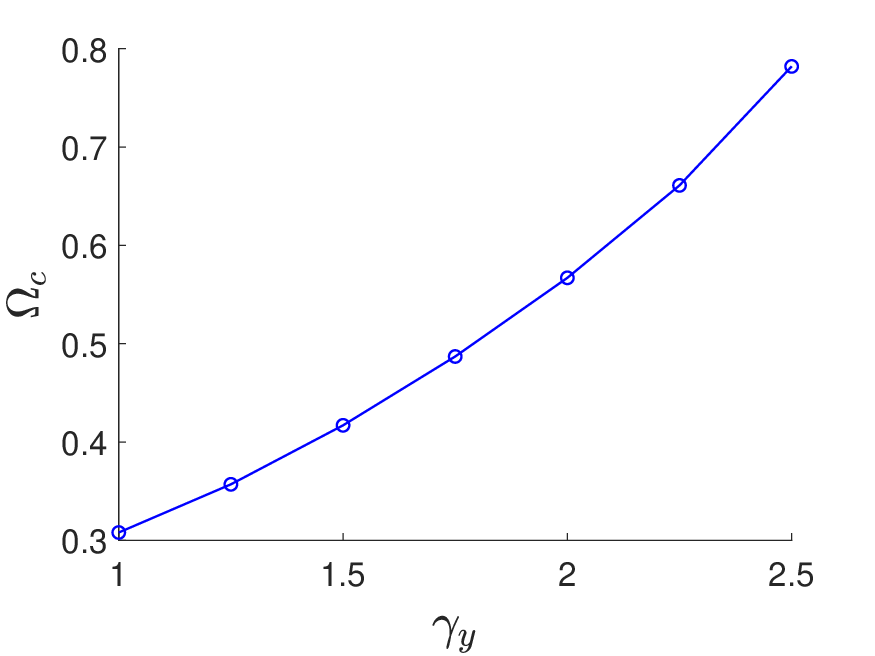,height=4.3cm,width=4.4cm,angle=0}\quad
\hspace{-0.4cm}
\psfig{figure=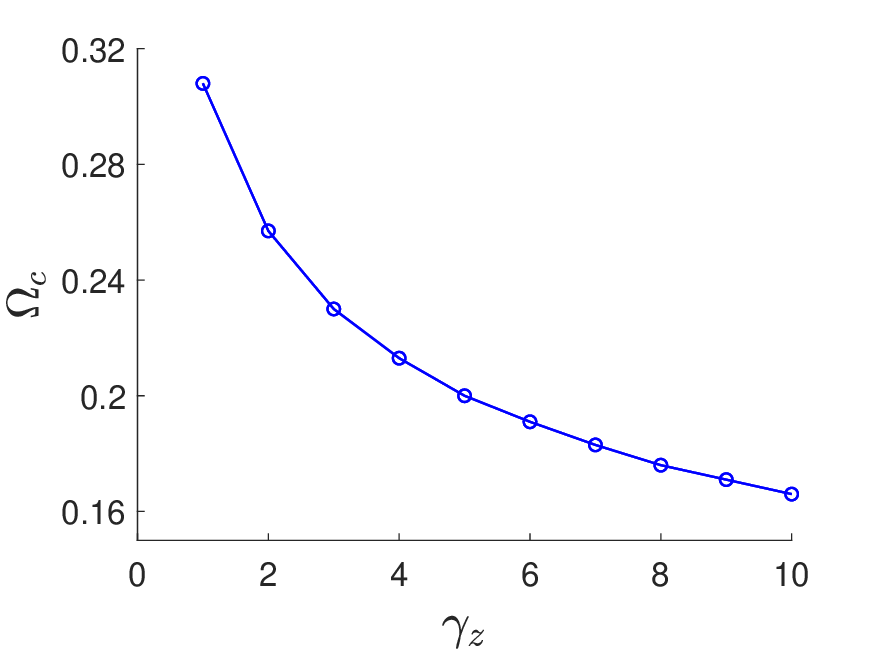,height=4.3cm,width=4.4cm,angle=0}\quad
\hspace{-0.4cm}
\psfig{figure=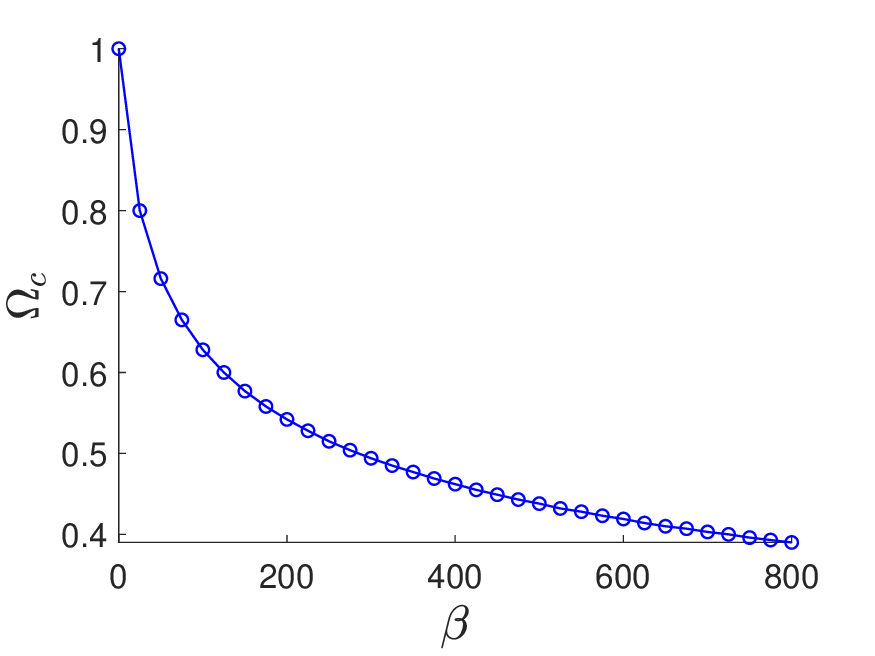,height=4.3cm,width=4.4cm,angle=0}}
\caption{The variations of $\Omega_{c}$ for \textbf{Case I-III} in \textbf{Example \ref{exp_critical1}}. }\label{critical1} 
\end{figure}

For \textbf{Case I-II}, the computation domain is chosen as $\mathcal{D}_{L\bm{\xi}}=[-16,16]^{3}$ with grid number $N=256^{3}$.
Fig. \ref{critical1}  depicts  the variations of critical rotational frequency $\Omega_{c}$ with respect to the trapping frequencies $\gamma_{y}$,
$\gamma_{z}$ and  local interaction strength $\beta$ respectively.
From this figure, we can observe clearly that: (i) The critical rotational frequency $\Omega_{c}$ increases with  $\gamma_y$,
in contrast to its decreasing trend with $\gamma_z$.
(ii) The critical rotational frequency $\Omega_{c}$ decreases as $\beta$ increases, and it decreases faster for smaller $\beta$.
Moreover, based on extensive numerical results not shown here, it is reasonable to conjecture that
$\Omega_{c}$ tends to one when $\beta$ tends to zero.

\begin{exmp} \label{exp_critical2}(\textbf{The impact of dipolar potential})
Given an isotropic harmonic potential,  i.e., $V(\bx)=\frac{1}{2}|\bx|^2$, we investigate the variations of $\Omega_c$ in the following cases:
\vspace{0.08cm}\begin{itemize}
\item []\textbf{Case I}: Let $\textbf{n}=(0, 0, 1)^T$,
take $\beta=200,300,400$ and vary $\lambda$;\vspace{0.05cm}
\item []\textbf{Case II}: Same as \textbf{Case I} except for $\textbf{n}=(1, 0, 0)^T$;\vspace{0.05cm}
\item []\textbf{Case III}: Let $\beta=400$,  $\textbf{n}=(\sin\theta,0,\cos\theta)^T$ with varying $\theta \in [0, \pi]$, 
    and take $\lambda=0,\pm50,\pm 100,\pm 200$ respectively.
\end{itemize}\vspace{0.08cm}
\end{exmp}

\begin{figure}[h!]
\centerline{
\psfig{figure=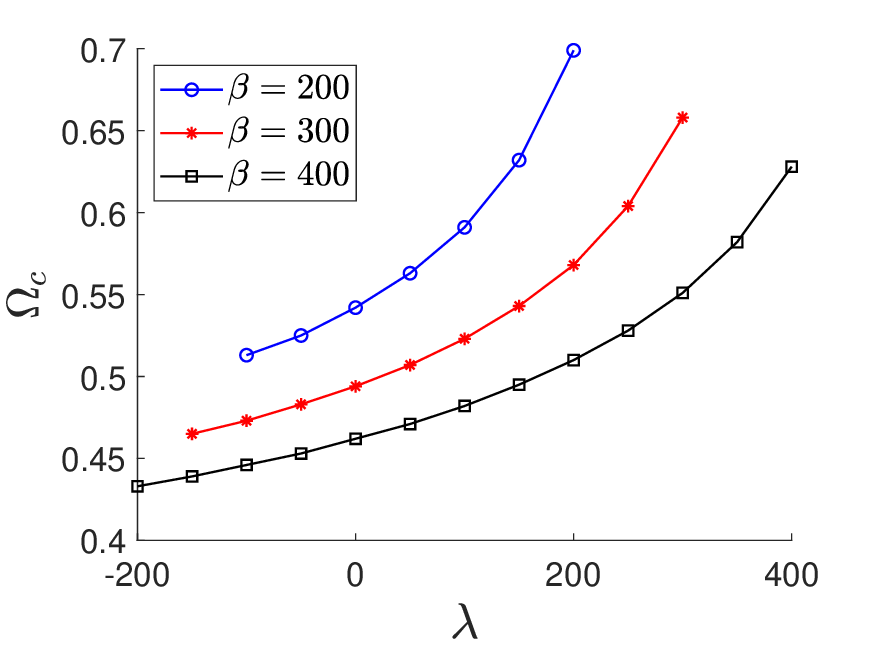,height=4.3cm,width=4.4cm,angle=0}\quad
\hspace{-0.43cm}
\psfig{figure=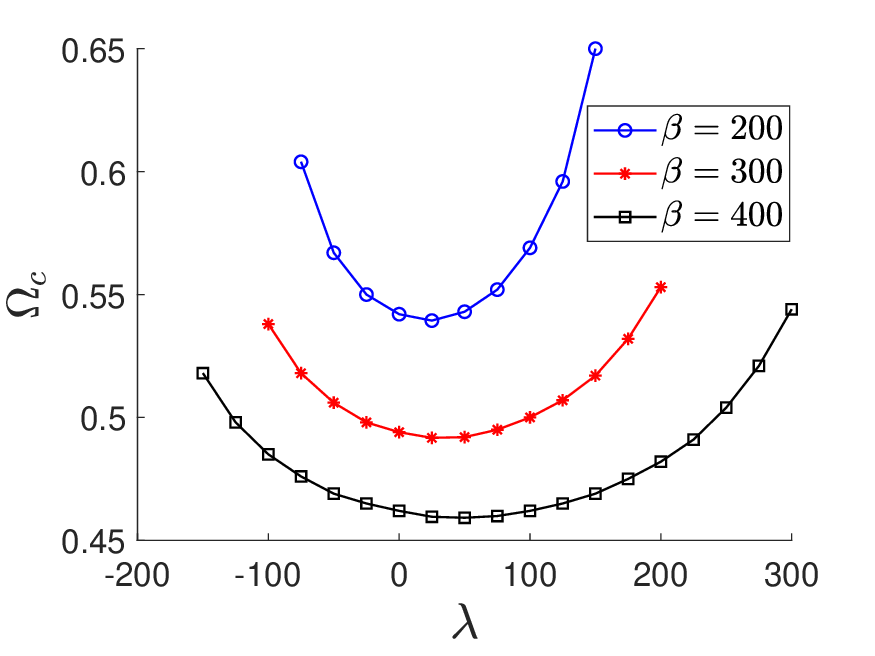,height=4.3cm,width=4.4cm,angle=0}\quad
\hspace{-0.43cm}
\psfig{figure=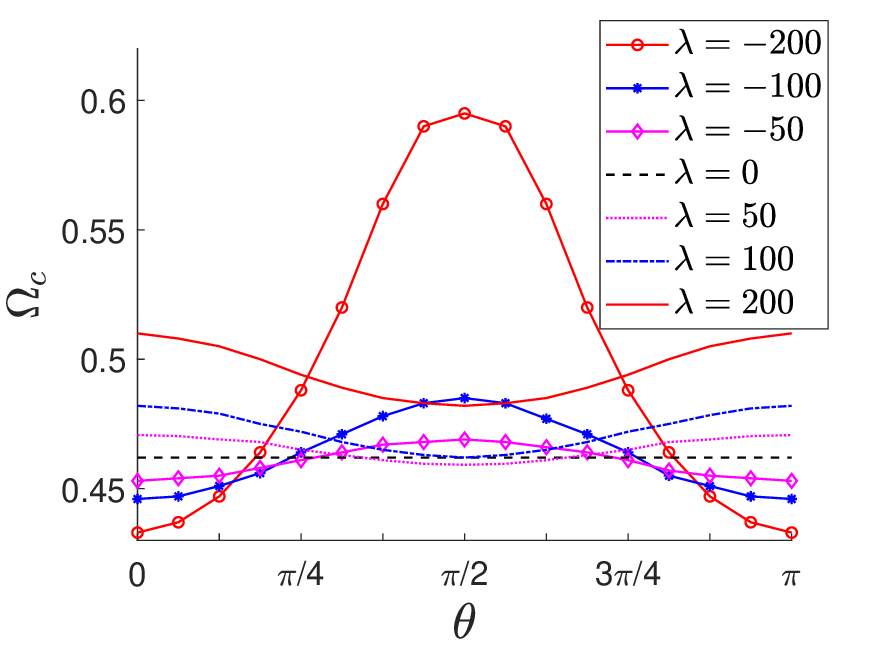,height=4.3cm,width=4.4cm,angle=0}}
\caption{The variations of $\Omega_{c}$ with $\lambda$ for $\bm{n}=(0, 0, 1)^T$, $(1, 0, 0)^T$ and that with $\theta$ in $\bm{n}=(\sin\theta,0,\cos\theta)^T$(from left to right) of \textbf{Example \ref{exp_critical2}}. }\label{critical_lambdaTheta} %
\end{figure}

Fig. \ref{critical_lambdaTheta} depicts the variations of critical rotational frequency $\Omega_{c}$  with respect to the dipolar potential strength $\lambda$
in \textbf{Case I-II}
and the angle $\theta$ in \textbf{Case III}.
From this figure, it is clear that $\Omega_{c}$ increases as $\lambda$ increases when the dipole orientation is $\bm{n}=(0,0,1)^T$.
For orientation $\bm{n}=(1,0,0)^T$,  $\Omega_{c}$ increases as $|\lambda|$ increases when $\lambda <0$. While, for positive $\lambda$,
it first decreases and then increases as $\lambda$ increases.
From the numerical results not shown here for brevity, given an isotropic harmonic potential $V(\bx)$, we conclude that:
(i) When the dipole orientation $\bm{n}$ lies in the $x$-$y$ plane, the variation behavior of $\Omega_c$ is the same as that of the case $\bm{n}=(1,0,0)^T$. (ii) There is no dependence of $\Omega_c$ on the azimuthal angle
of the orientation $\bm{n}$.


In the third picture of Fig. \ref{critical_lambdaTheta}, $\theta$ denotes the angle between the dipole orientation $\bm{n}$ and the rotation axis ($z$-axis).
For  positive $\lambda$, $\Omega_{c}$ decreases as the dipole orientation rotates from parallel ($\theta=0, \pi$) to  perpendicular ($\theta=\pi/2$) alignment relative  to the rotation axis,
while for negative $\lambda$, it increases.
Moreover, the magnitude of this variation in $\Omega_{c}$ is smaller for a given positive $\lambda=\alpha$  than for  its negative counterpart $\lambda=-\alpha$.

\subsection{Impacts on energies and chemical potential}
In this subsection, we study the variation trends of chemical potential $\mu$ and energies $E$, $E_{\rm kin}$, $E_{\rm pot}$, $E_{\rm int}$, $E_{\rm dip}$,
$E_{\rm rot}:=-\Omega\int_{{\mathbb R}^3} \bar{\psi} L_z\psi {\rm d} \bx$
with respect to different physical parameters, including trapping frequencies of the harmonic potential \eqref{Vpoten},
local interaction strength $\beta$, dipolar potential strength $\lambda$ and rotational frequency $\Omega$.
We choose computation domain $\mathcal{D}_{L\bm{\xi}}=[-20,20]^{3}$ and take grid number $N=256^{3}$, if not stated otherwise.

\begin{exmp} \label{Emu_1}(\textbf{The impact of external potential})
Here, we investigate the dependence of the energies and chemical potential on the external potential. To this end, we choose the harmonic potential with $\gamma_y=1$, $\beta=600$, $\Omega=0.6$, $\lambda=-200$,  $\textbf{n}=(0, 0, 1)^T$ and consider the following cases:
\vspace{0.08cm}\begin{itemize}
\item []\textbf{Case I}: Let $\gamma_z=1$  and vary $\gamma_x$;\vspace{0.05cm}
\item []\textbf{Case II}: Let $\gamma_x=1$  and vary $\gamma_z$.
\end{itemize}\vspace{0.08cm}
\end{exmp}

\begin{figure}[h!]
\centerline{
\psfig{figure=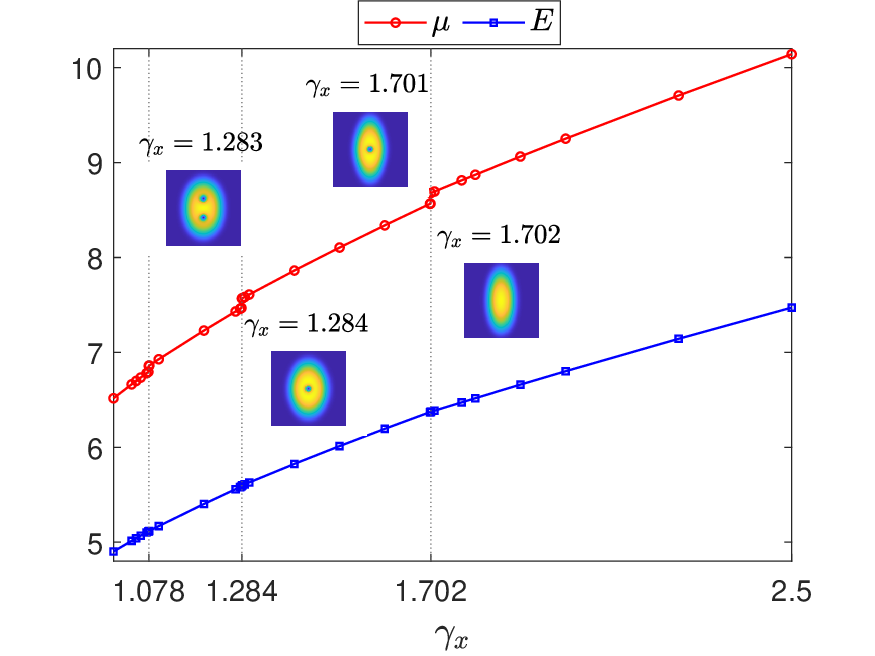,height=5.4cm,width=7.3cm,angle=0}\quad
\hspace{-0.9cm}
\psfig{figure=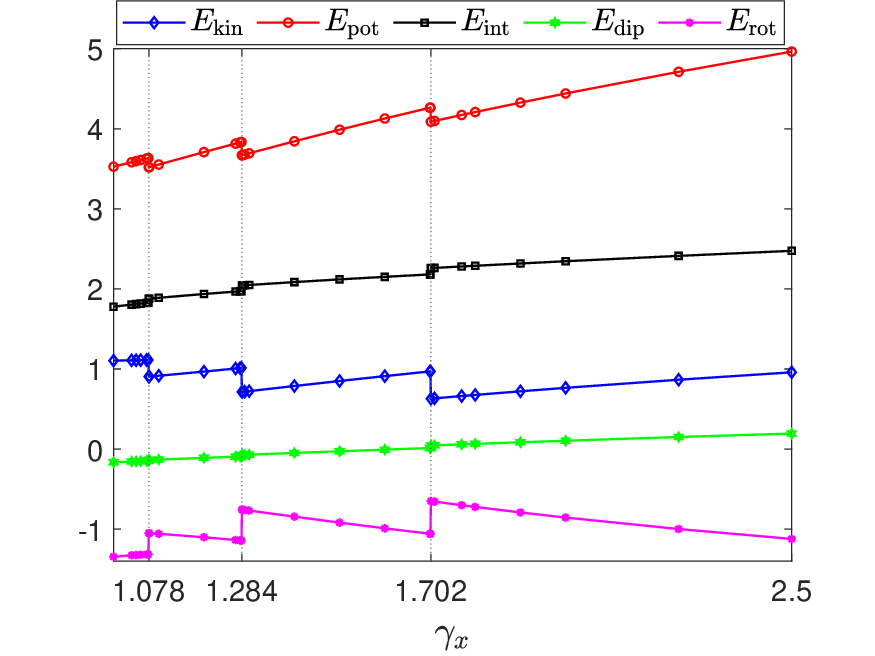,height=5.4cm,width=7.3cm,angle=0}}
\vspace{0.3cm}

\centerline{
\psfig{figure=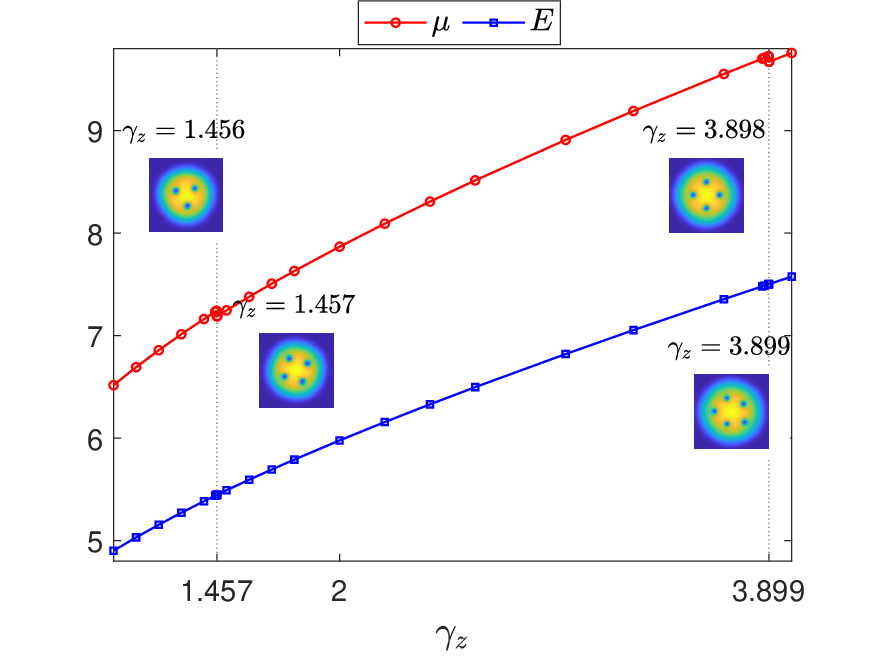,height=5.4cm,width=7.3cm,angle=0}\quad
\hspace{-0.9cm}
\psfig{figure=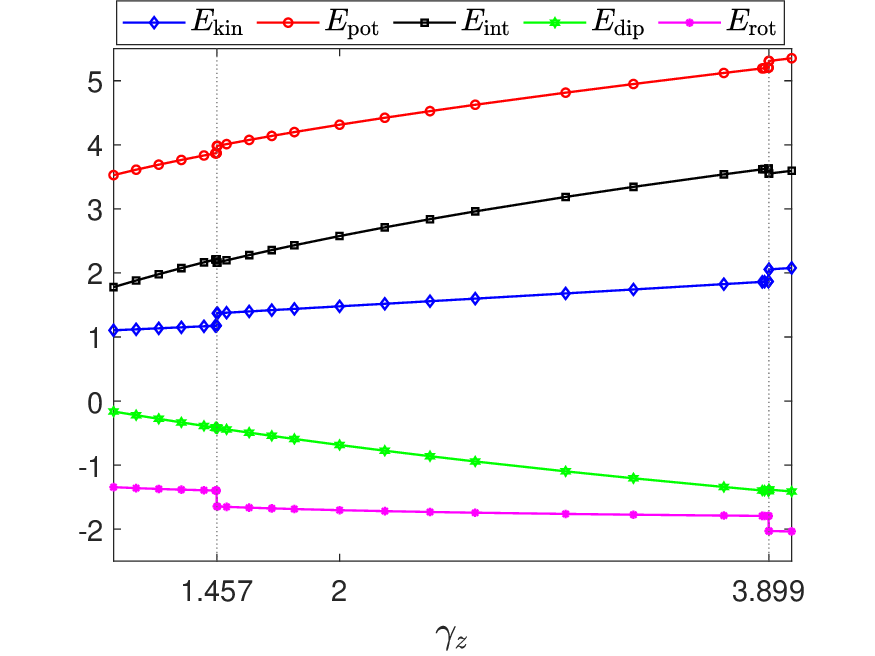,height=5.4cm,width=7.3cm,angle=0}}
\caption{The variations of energies and $\mu$ for \textbf{Case I-II} (from top to bottom) in \textbf{Example \ref{Emu_1}}.}\label{Emu_fig1}
\end{figure}

Fig. \ref{Emu_fig1} depicts the variations of energies and chemical potential
in \textbf{Case I-II}, from which we can see  that:
(i) Only the energies $E$ and $E_{\rm dip}$ vary continuously,
whereas chemical potential $\mu$ and $E_{\rm kin}$, $E_{\rm pot}$, $E_{\rm int}$, $E_{\rm rot}$ exhibit `jumps'
at the phase transition point where  the vortex number changes.
(ii) The energy $E$  exhibits a monotonic increase  with respect to both $\gamma_x$ and $\gamma_z$.
(iii) The energy $E_{\rm dip}$ exhibits a monotonic increase (decrease)  with respect to  $\gamma_x$ ($\gamma_z$).

\begin{figure}[h!]
\centerline{
\psfig{figure=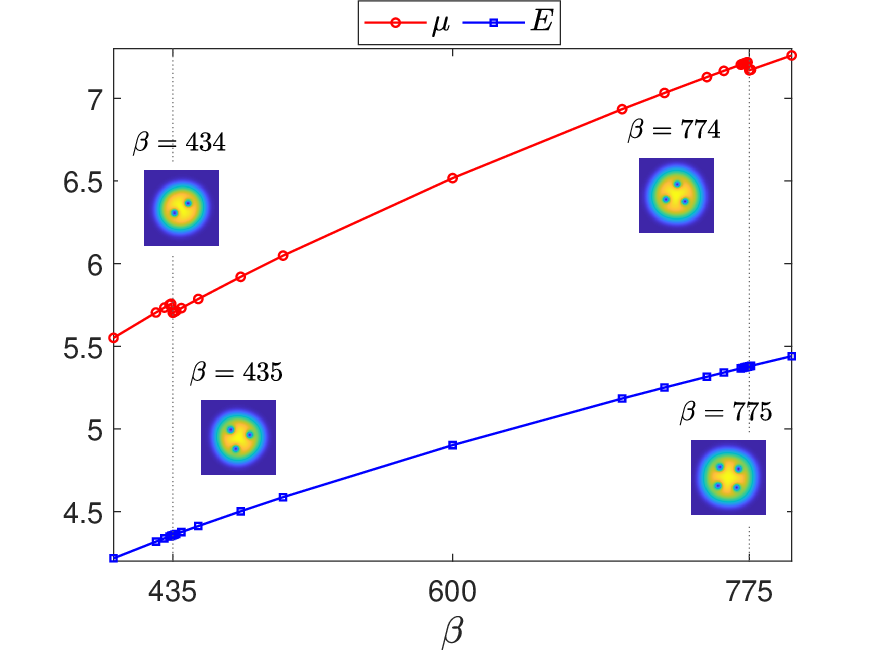,height=5.4cm,width=7.3cm,angle=0}\quad
\hspace{-0.9cm}
\psfig{figure=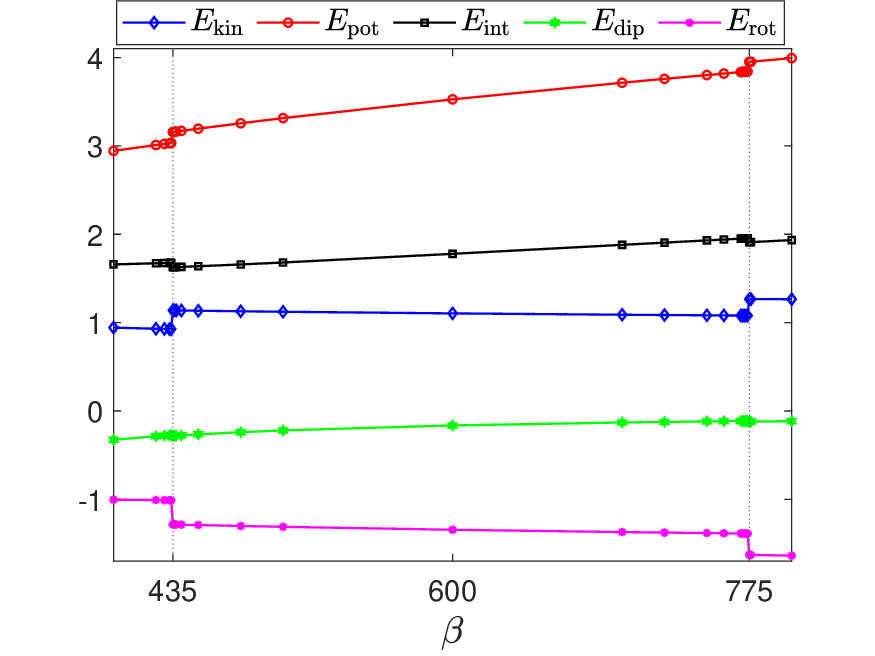,height=5.4cm,width=7.3cm,angle=0}}
\vspace{0.3cm}

\centerline{
\psfig{figure=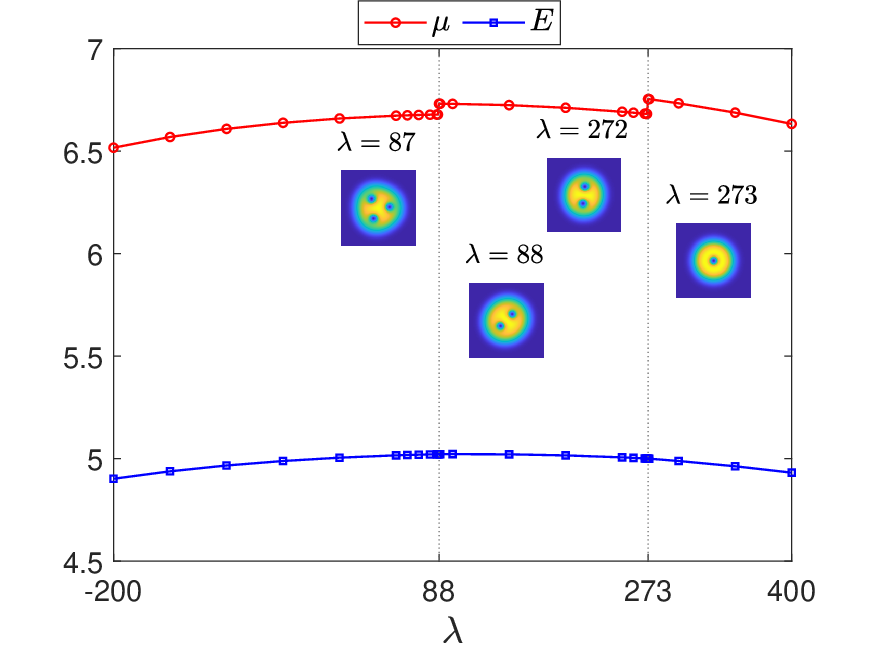,height=5.4cm,width=7.3cm,angle=0}\quad
\hspace{-0.9cm}
\psfig{figure=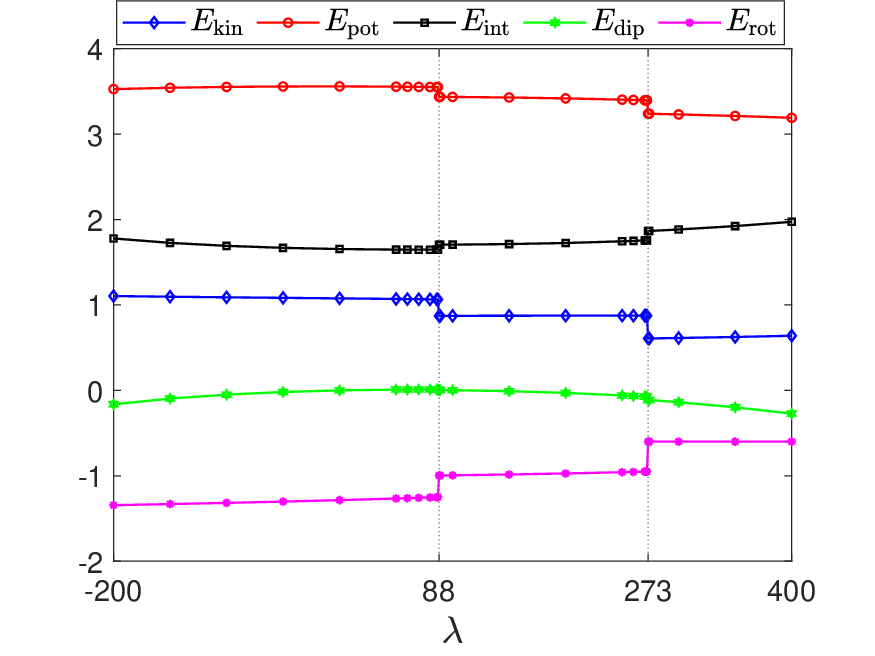,height=5.4cm,width=7.3cm,angle=0}}
\vspace{0.3cm}

\centerline{
\psfig{figure=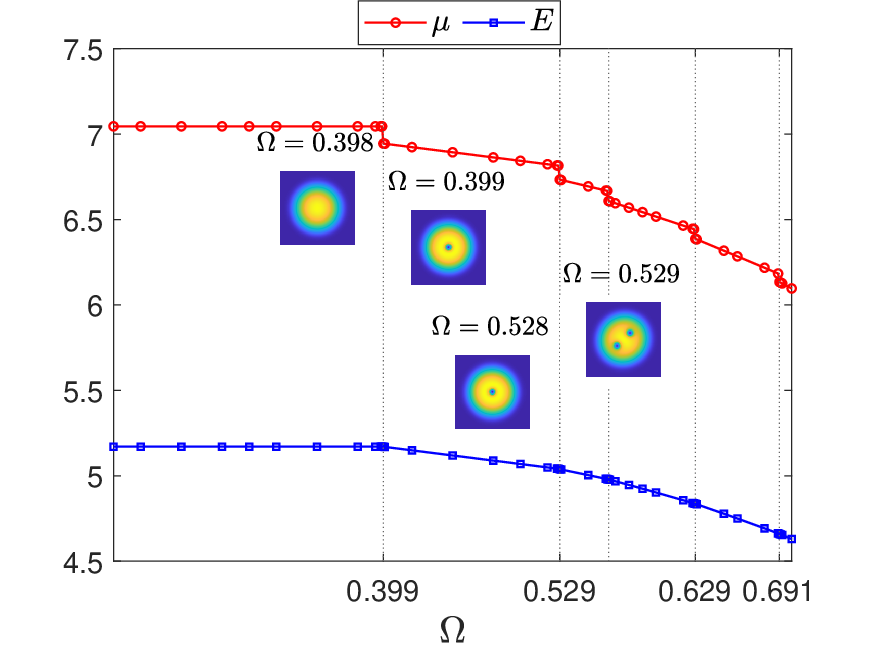,height=5.4cm,width=7.3cm,angle=0}\quad
\hspace{-0.9cm}
\psfig{figure=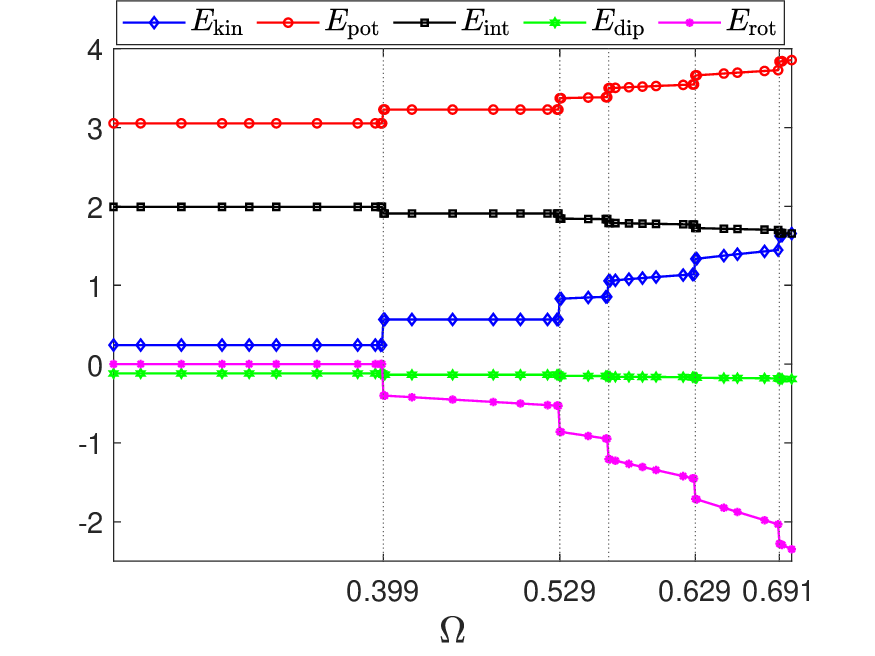,height=5.4cm,width=7.3cm,angle=0}}
\caption{The variations of energies and $\mu$ for \textbf{Case I-III} (from top to bottom) in \textbf{Example \ref{Emu_2}}.}\label{Emu_fig2}
\end{figure}
\begin{exmp} \label{Emu_2}(\textbf{The impacts of short/long-range interaction and rotational frequency})
Given the dipole orientation $\textbf{n}=(0, 0, 1)^T$ and an isotropic harmonic potential,  i.e., $V(\bx)=\frac{1}{2}|\bx|^2$, we investigate the variations of energies and chemical potential in the following cases:
\vspace{0.08cm}\begin{itemize}
\item []\textbf{Case I}: Let $\Omega=0.6$, $\lambda=-200$, and vary $\beta$;\vspace{0.05cm}
\item []\textbf{Case II}: Let $\Omega=0.6$, $\beta=600$, and vary $\lambda$;\vspace{0.05cm}
\item []\textbf{Case III}: Let $\beta=600$, $\lambda=-200$, and vary $\Omega$.
\end{itemize}\vspace{0.08cm}
\end{exmp}

We take computation domain $\mathcal{D}_{L\bm{\xi}}=[-16,16]^{3}$  for all cases.
Fig. \ref{Emu_fig2} depicts the variations of energies and chemical potential, from which we can see clearly that:
(i) Only the energies $E$ and $E_{\rm dip}$ vary continuously, whereas chemical potential $\mu$ and
$E_{\rm kin}, E_{\rm pot}, E_{\rm int}, E_{\rm rot}$ display abrupt discontinuities in the presence of phase transition in the ground state.
(ii) Both  $E$ and  $E_{\rm dip}$ increase (decrease) monotonically with respect to  $\beta$ ($\Omega$).
(iii) Under an isotropic trap, when there is no vortex in the ground state, all the energies and chemical potential $\mu$
remain unchanged as $\Omega$ increases.

\subsection{Patterns of ground states in anisotropic traps}
Here we compute the ground states of anisotropic rotating dipolar BEC and present the corresponding ground state patterns.
In the following numerical tests, we take $\beta=5000$, $\Omega=0.85$, $\gamma_{y}=\gamma_{z}=1$ and choose computation domain $\mathcal{D}_{L\bm{\xi}}=[-14,14]^{3}$  with grid number $N=256^{3}$.

\begin{exmp}\label{ex_withoutDDI}
(\textbf{Non-dipolar BEC})
In the absence of dipolar potential, i.e., $\lambda=0$, we take $\gamma_{x}=1, 1.5, 2, 3$ respectively and compute the  ground states.
\end{exmp}

\begin{figure}[h!]
\centerline{
\psfig{figure=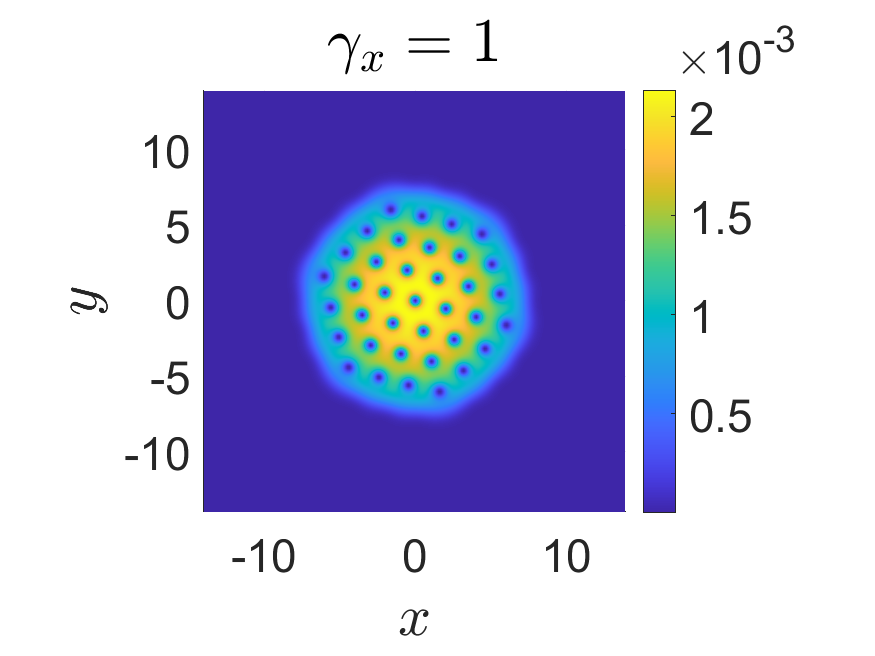,height=2.8cm,width=3.7cm,angle=0}\quad
\hspace{-0.8cm}
\psfig{figure=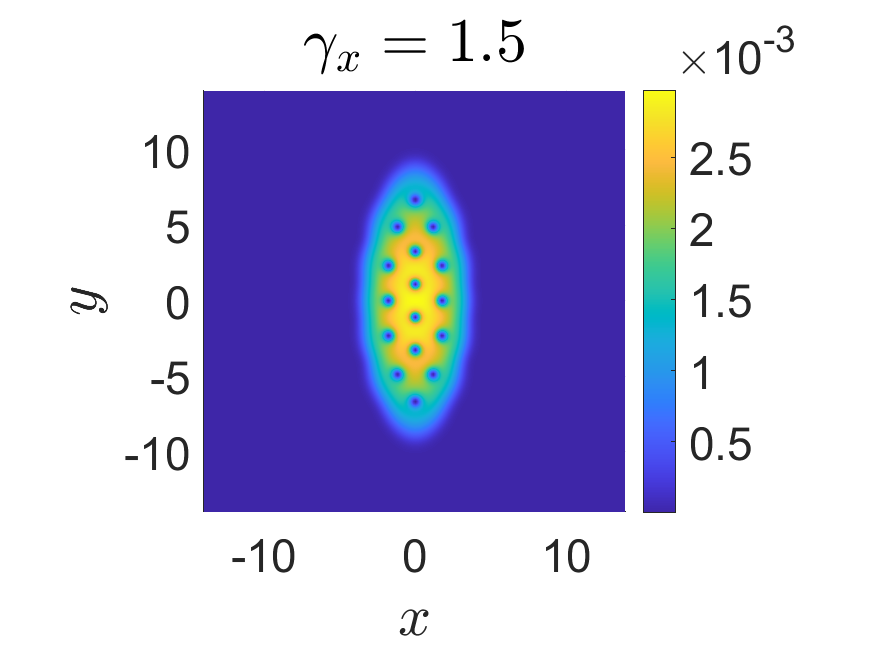,height=2.8cm,width=3.7cm,angle=0}\quad
\hspace{-0.8cm}
\psfig{figure=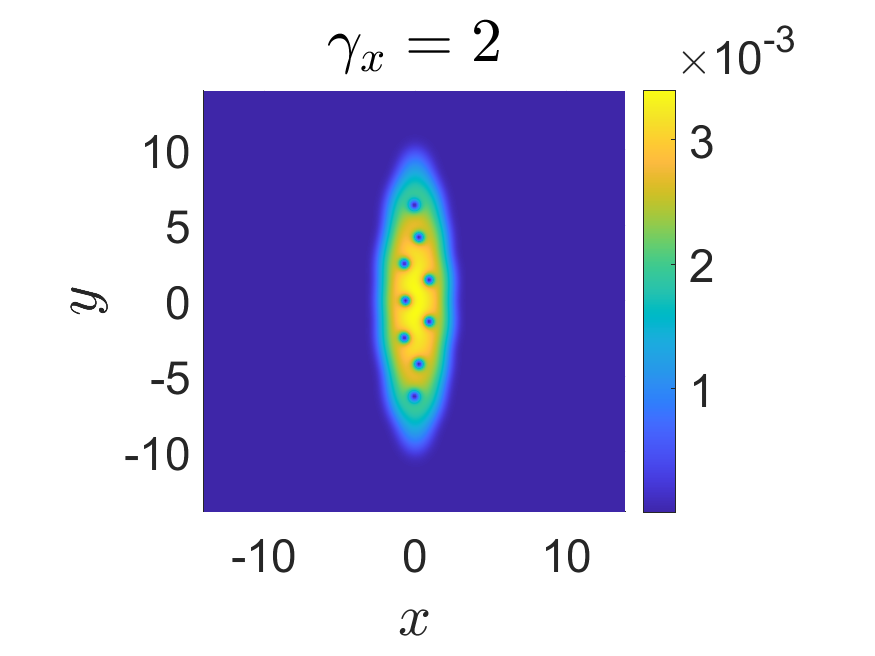,height=2.8cm,width=3.7cm,angle=0}\quad
\hspace{-0.8cm}
\psfig{figure=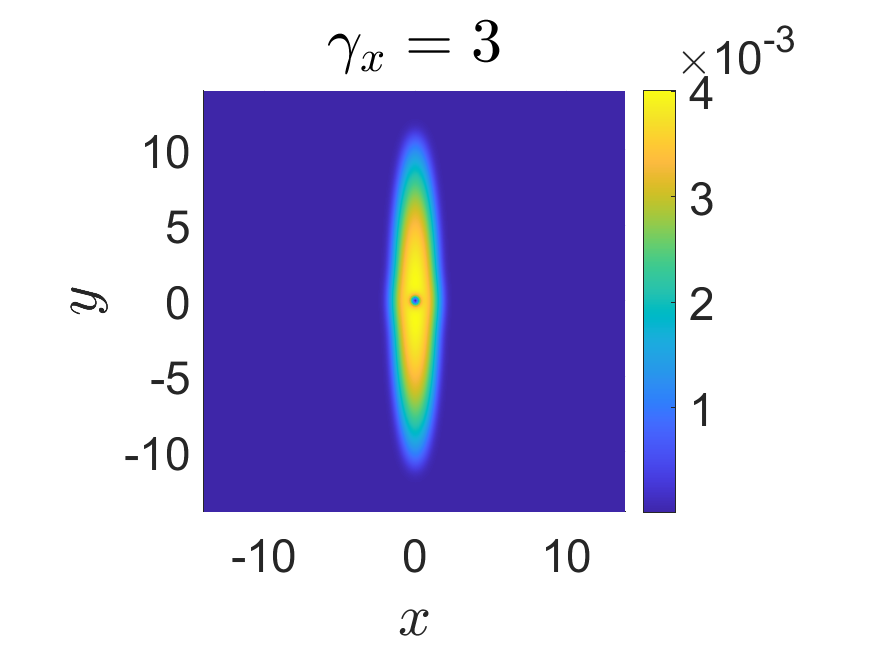,height=2.8cm,width=3.7cm,angle=0}}
\vspace{-0.1cm}
\centerline{
\psfig{figure=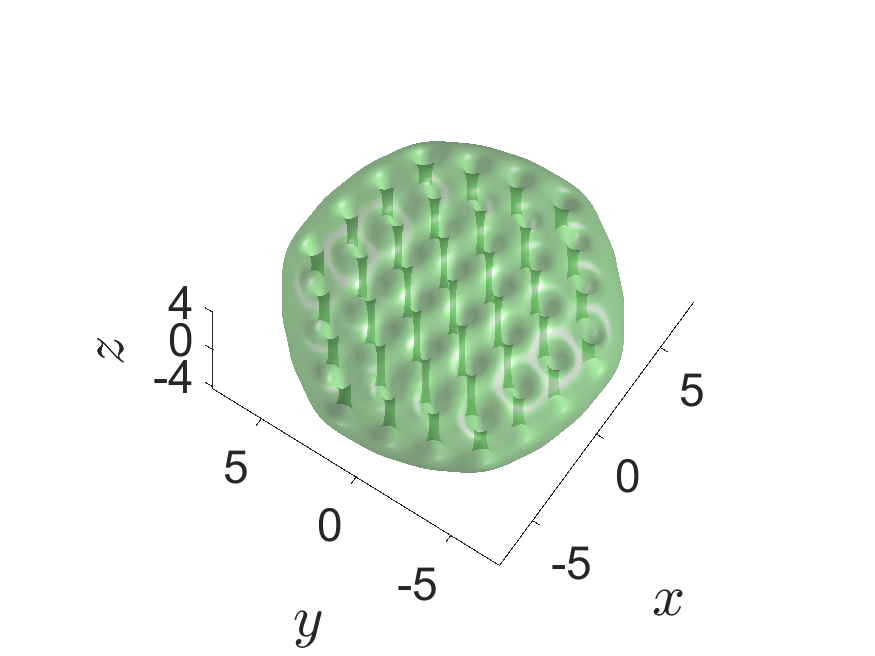,height=2.8cm,width=3.7cm,angle=0}\quad
\hspace{-0.8cm}
\psfig{figure=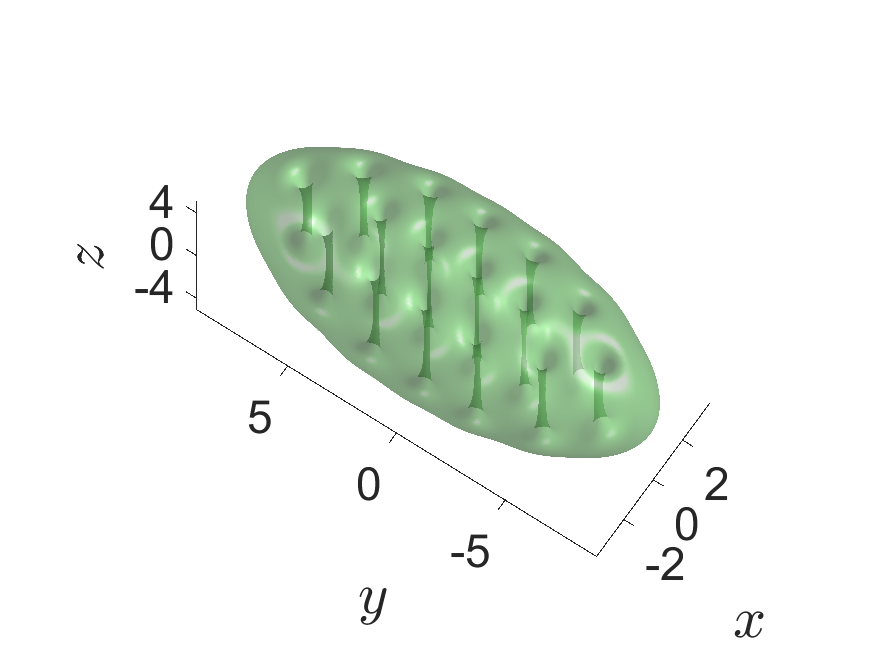,height=2.8cm,width=3.7cm,angle=0}\quad
\hspace{-0.8cm}
\psfig{figure=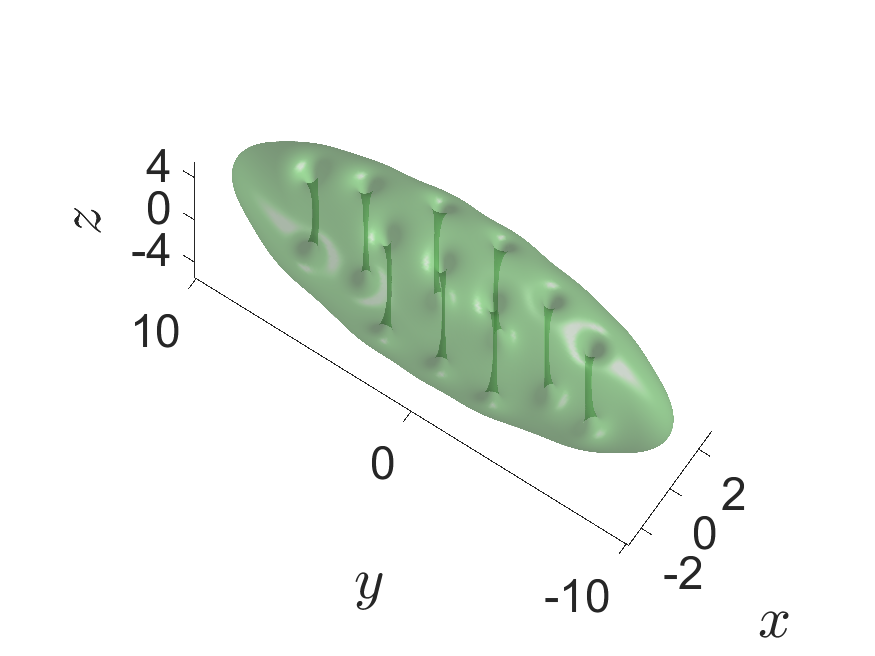,height=2.8cm,width=3.7cm,angle=0}\quad
\hspace{-0.8cm}
\psfig{figure=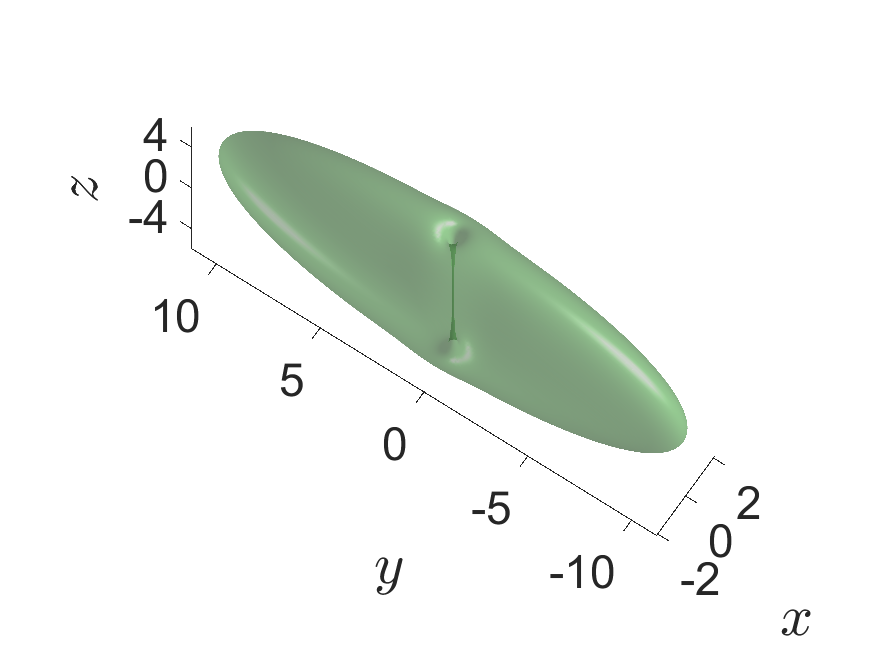,height=2.8cm,width=3.7cm,angle=0}}
\caption{Contour plots of $|\phi_g(x, y, z=0)|^2$ (top row) and  isosurface plots of $|\phi_g(\bx)|^2=3\times10^{-4}$ (bottom row)  in \textbf{Example \ref{ex_withoutDDI}}.}\label{vortice_noDDI}
\end{figure}

Fig. \ref{vortice_noDDI} depicts the contour plots of $|\phi_g(x, y, z=0)|^2$ and corresponding isosurface plots of $|\phi_g(\bx)|^2=3\times10^{-4}$ for different $\gamma_{x}$, from which we can see that:
(i) The vortex lines are parallel to the rotation axis.
(ii) Under an isotropic harmonic potential, the vortex lattice exhibits a triangular arrangement characteristic of the
Abrikosov lattice (cf. the first column in Fig. \ref{vortice_noDDI}).
(iii) As $\gamma_{x}$ increases, i.e., the anisotropy strength becomes stronger, the density profile gradually narrows in the $x$-direction and the number of vortices decreases.

\begin{exmp}\label{ex_withDDI}
(\textbf{Dipolar BEC}) Taking $\gamma_{x}=1, 1.5, 2, 3$ respectively, we compute the  ground states in the following cases:
\vspace{0.08cm}\begin{itemize}
\item []\textbf{Case I}: Let $\lambda=-2000$ and different dipole orientations:
\threeitems{[]$\qquad$}{$(a)$  $\textbf{n}=(1, 0, 0)^T$;}{$(b)$  $\textbf{n}=(0, 0, 1)^T$;}
\vspace{0.05cm}
\item []\textbf{Case II}: Let $\lambda=2000$ and different dipole orientations:
\threeitems{[]$\qquad$}{$(c)$  $\textbf{n}=(1, 0, 0)^T$;}{$(d)$  $\textbf{n}=(1/2, 1/2, \sqrt{2}/2)^T$.}
\end{itemize}\vspace{0.08cm}
\end{exmp}

\begin{figure}[h!]
\centerline{
\psfig{figure=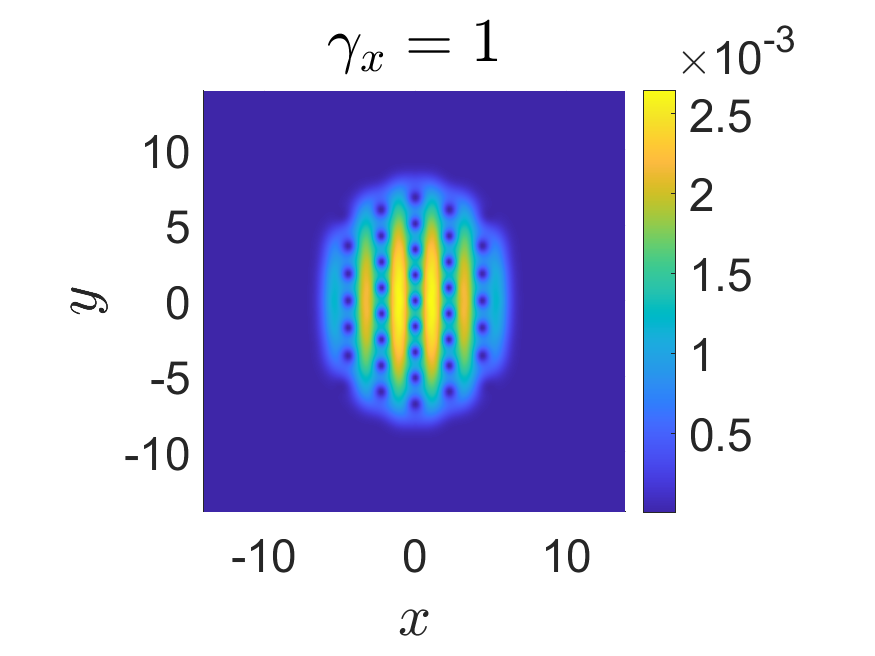,height=2.8cm,width=3.7cm,angle=0}\quad
\hspace{-0.8cm}
\psfig{figure=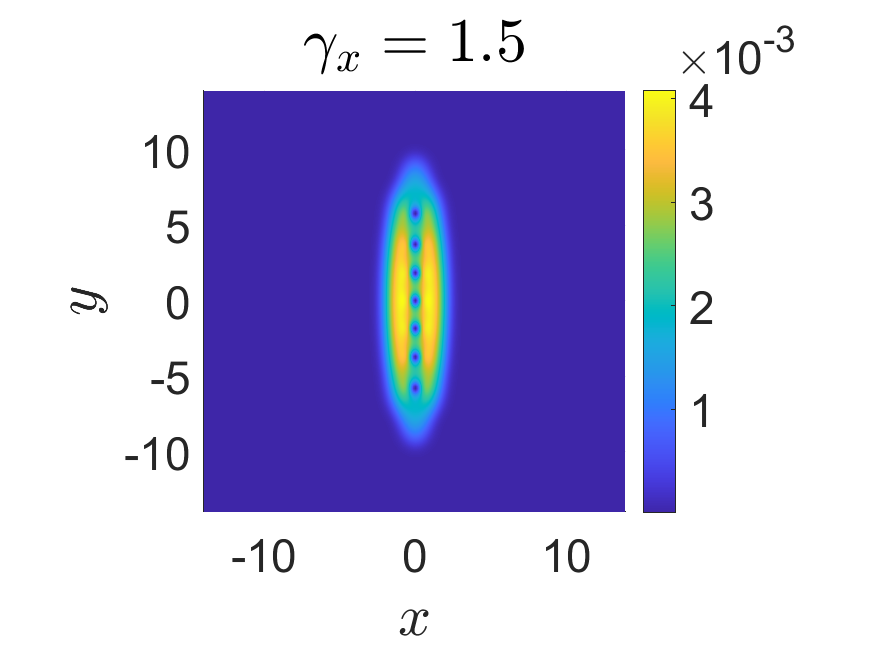,height=2.8cm,width=3.7cm,angle=0}\quad
\hspace{-0.8cm}
\psfig{figure=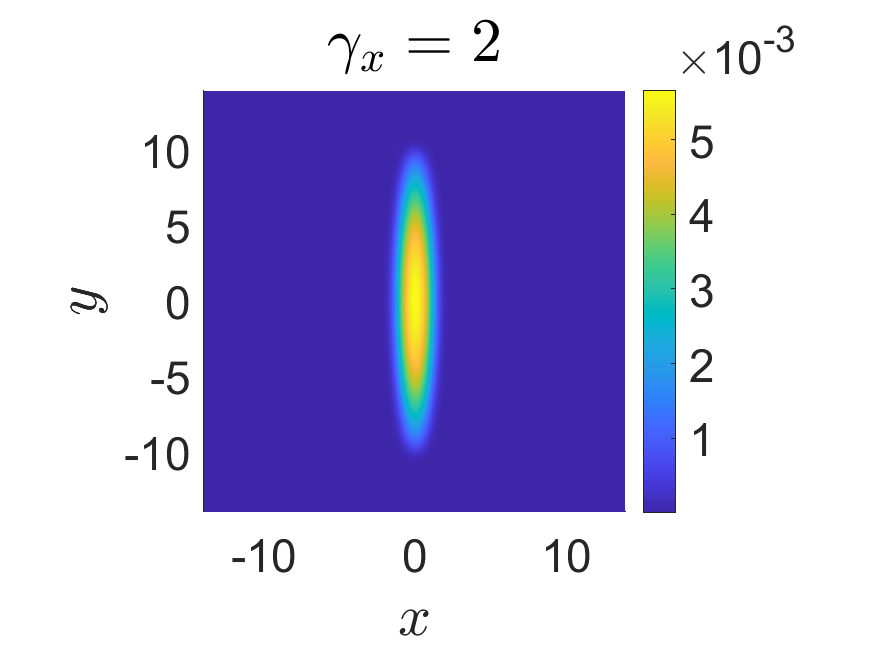,height=2.8cm,width=3.7cm,angle=0}\quad
\hspace{-0.8cm}
\psfig{figure=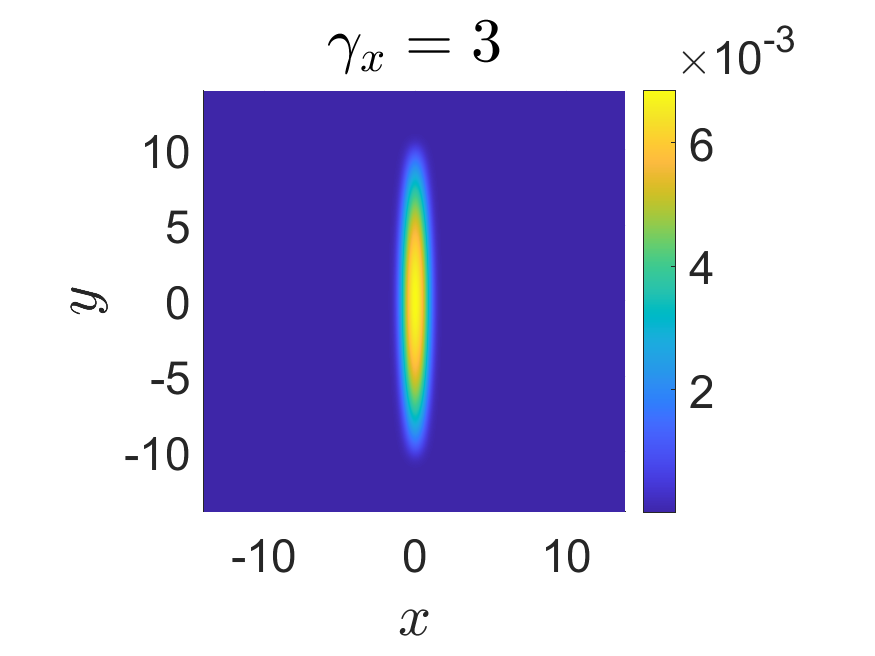,height=2.8cm,width=3.7cm,angle=0}}

\centerline{
\psfig{figure=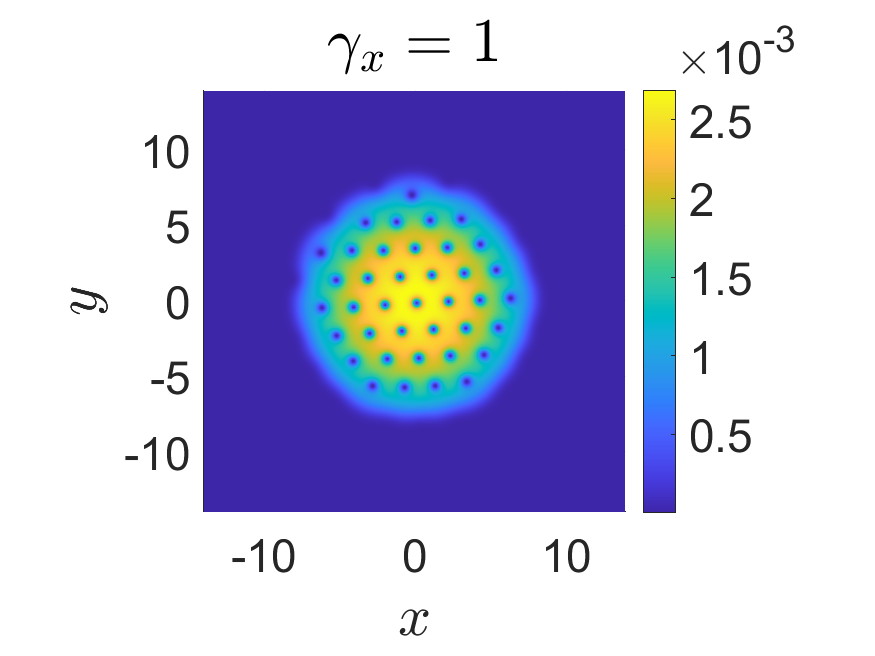,height=2.8cm,width=3.7cm,angle=0}\quad
\hspace{-0.8cm}
\psfig{figure=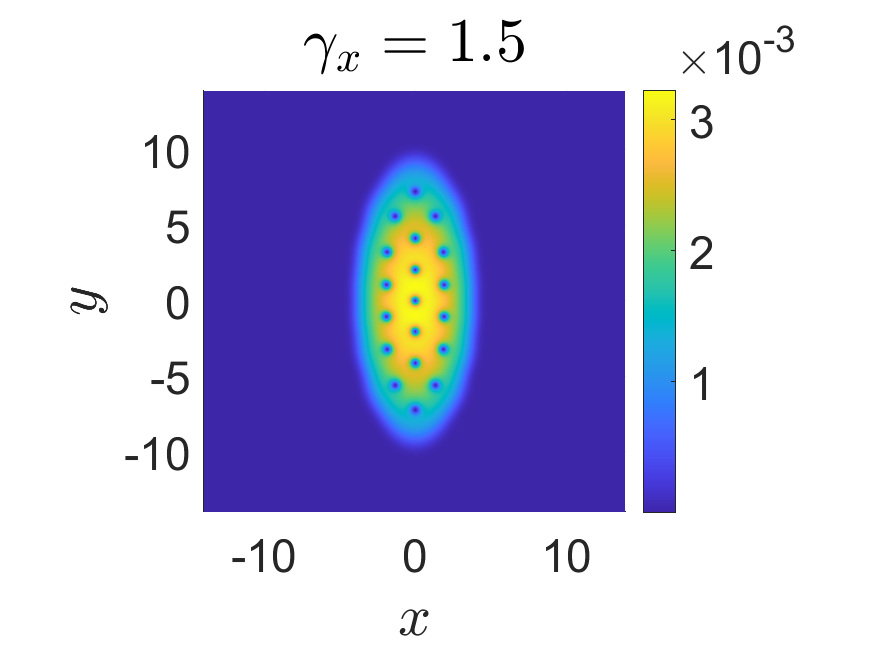,height=2.8cm,width=3.7cm,angle=0}\quad
\hspace{-0.8cm}
\psfig{figure=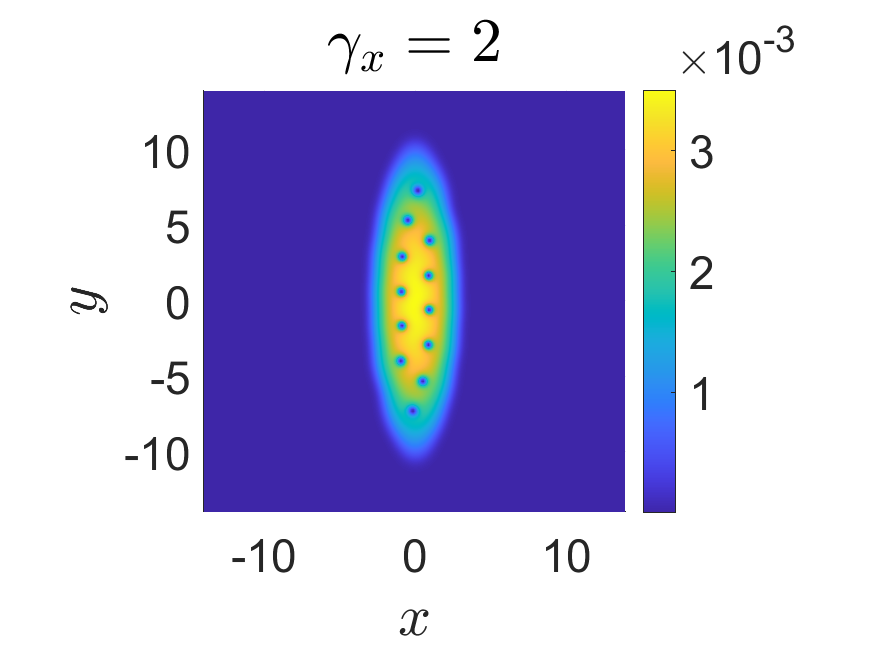,height=2.8cm,width=3.7cm,angle=0}\quad
\hspace{-0.8cm}
\psfig{figure=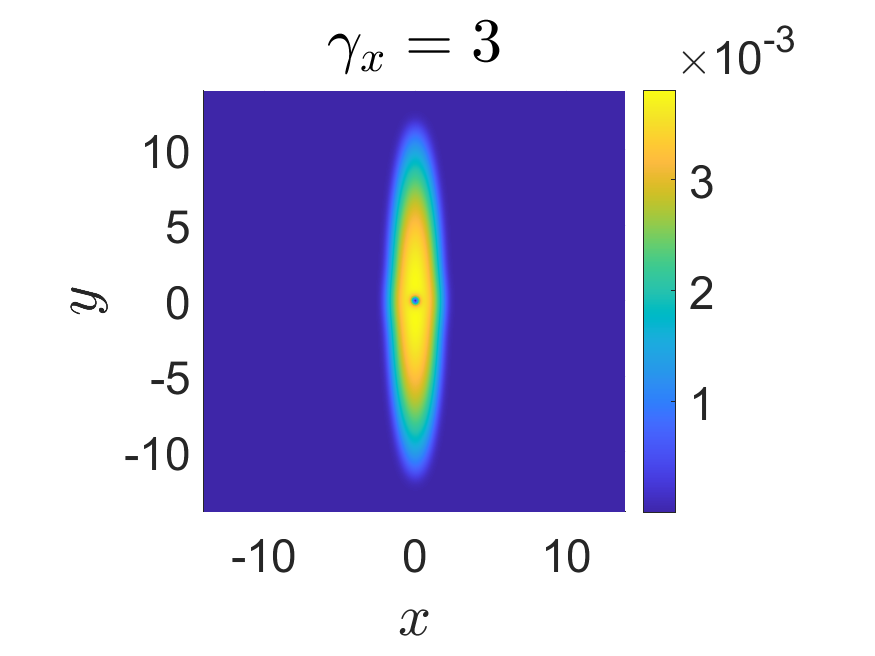,height=2.8cm,width=3.7cm,angle=0}}

\centerline{
\psfig{figure=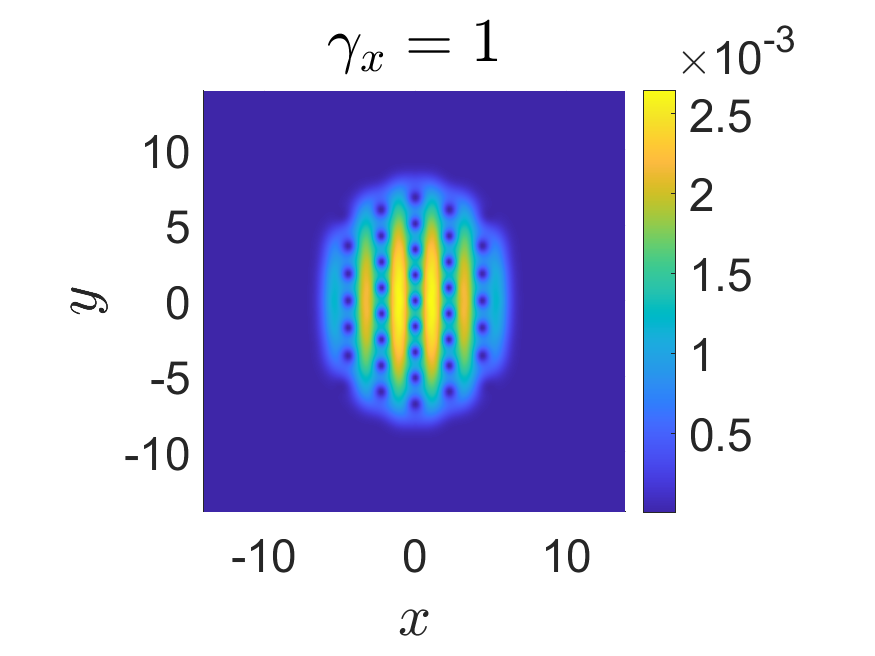,height=2.8cm,width=3.7cm,angle=0}\quad
\hspace{-0.8cm}
\psfig{figure=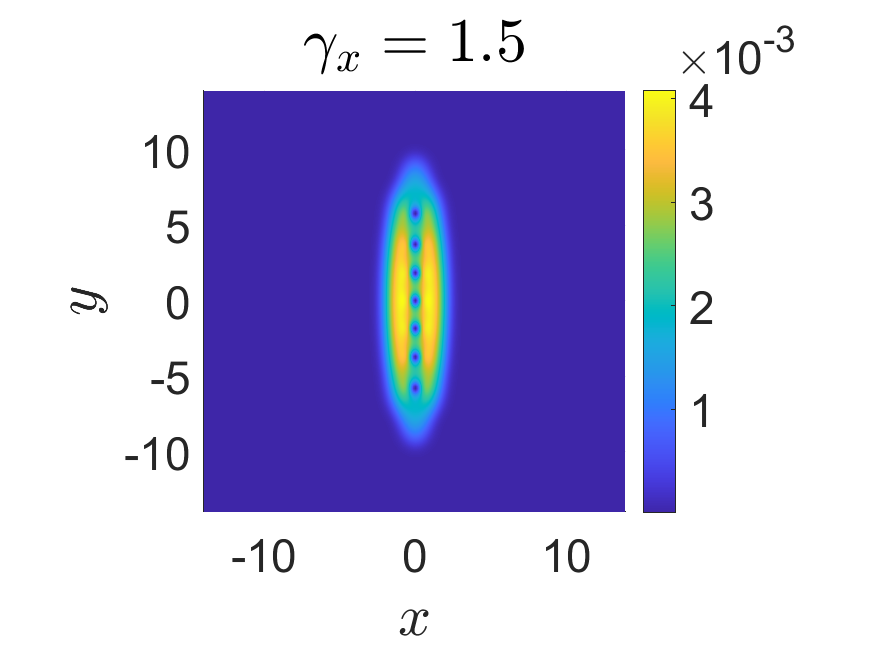,height=2.8cm,width=3.7cm,angle=0}\quad
\hspace{-0.8cm}
\psfig{figure=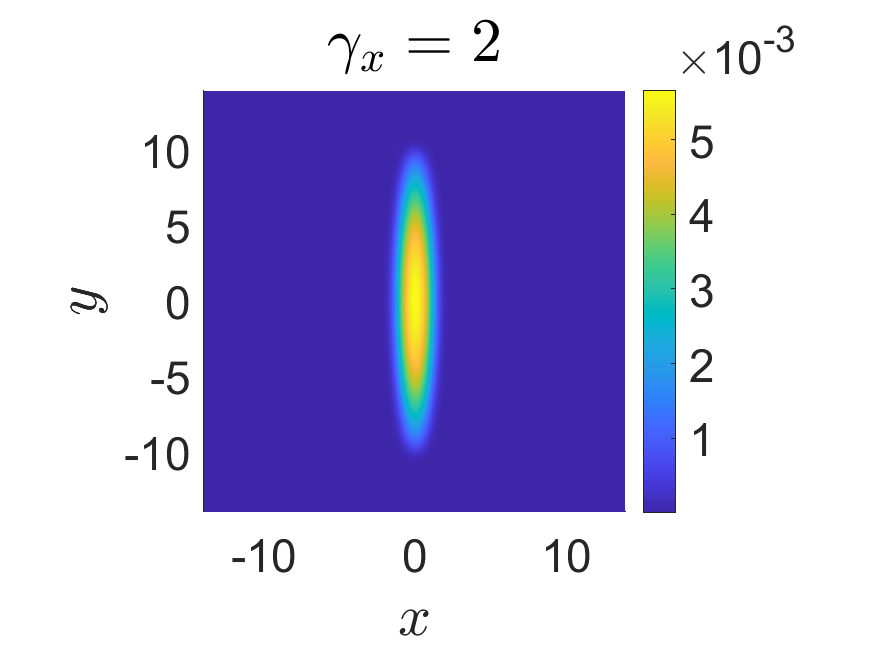,height=2.8cm,width=3.7cm,angle=0}\quad
\hspace{-0.8cm}
\psfig{figure=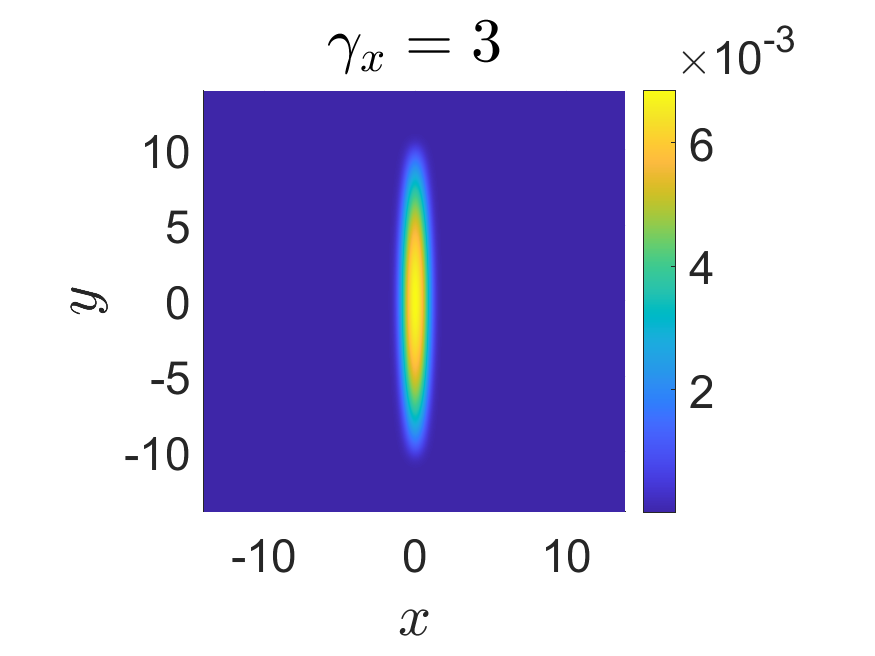,height=2.8cm,width=3.7cm,angle=0}}

\centerline{
\psfig{figure=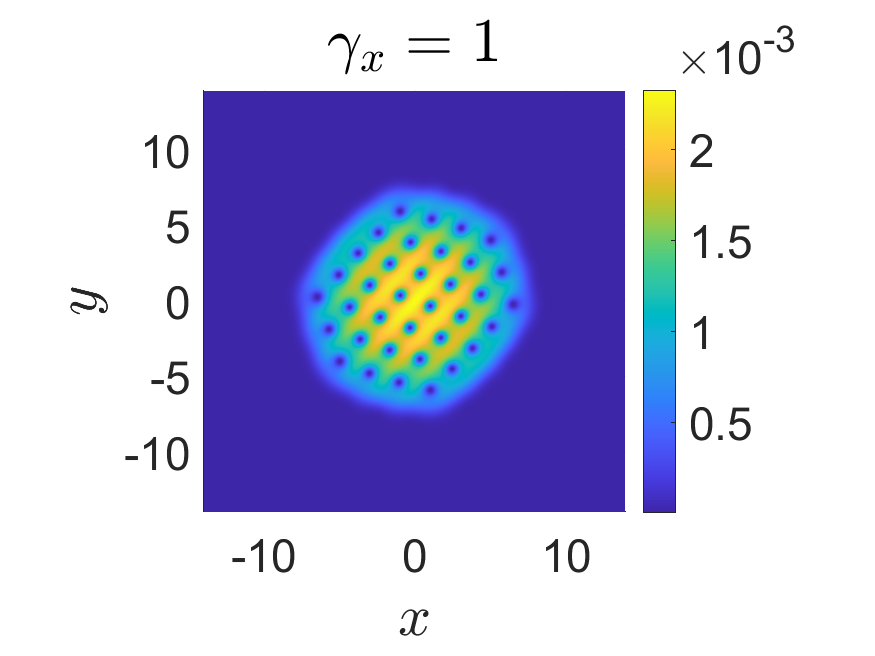,height=2.8cm,width=3.7cm,angle=0}\quad
\hspace{-0.8cm}
\psfig{figure=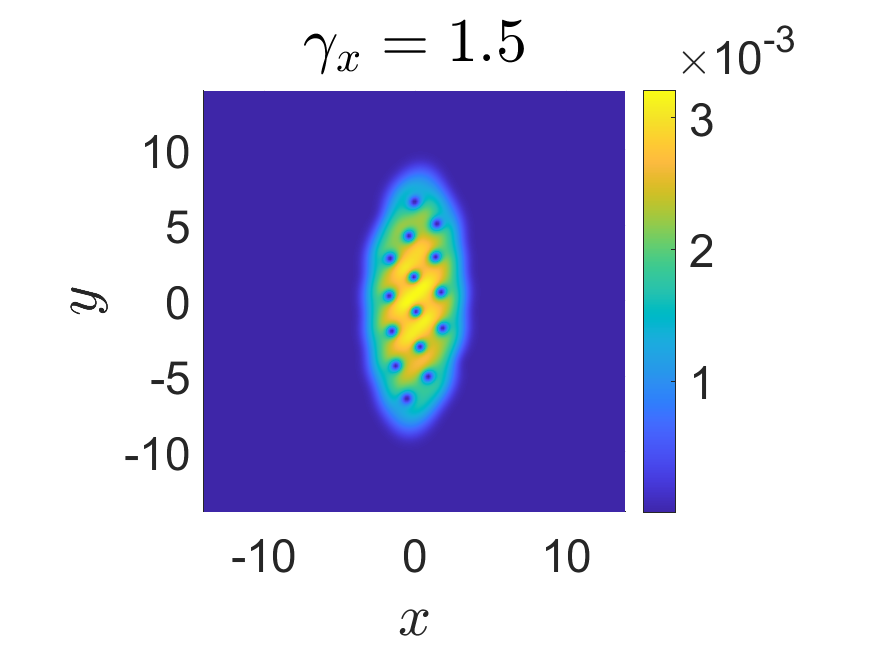,height=2.8cm,width=3.7cm,angle=0}\quad
\hspace{-0.8cm}
\psfig{figure=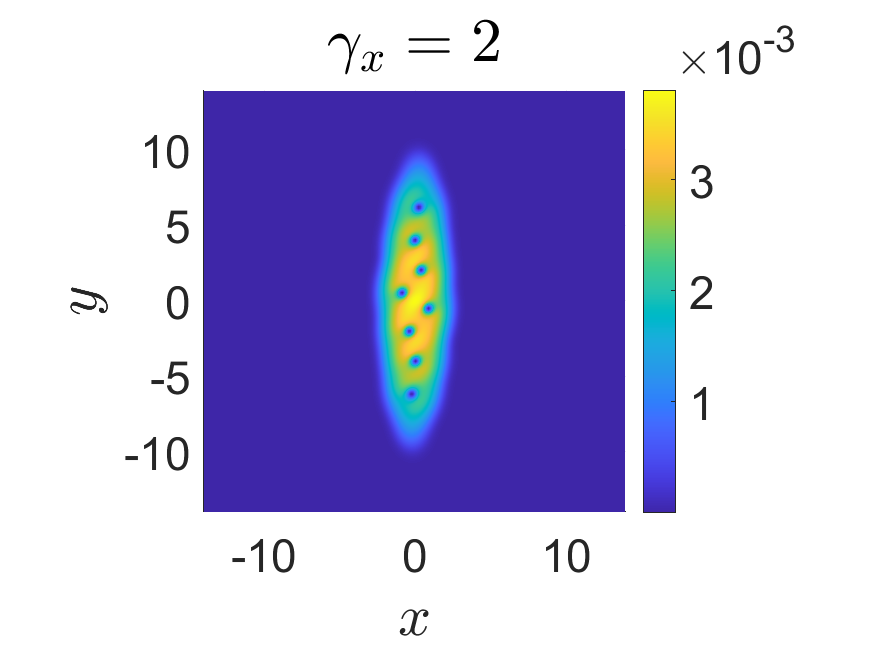,height=2.8cm,width=3.7cm,angle=0}\quad
\hspace{-0.8cm}
\psfig{figure=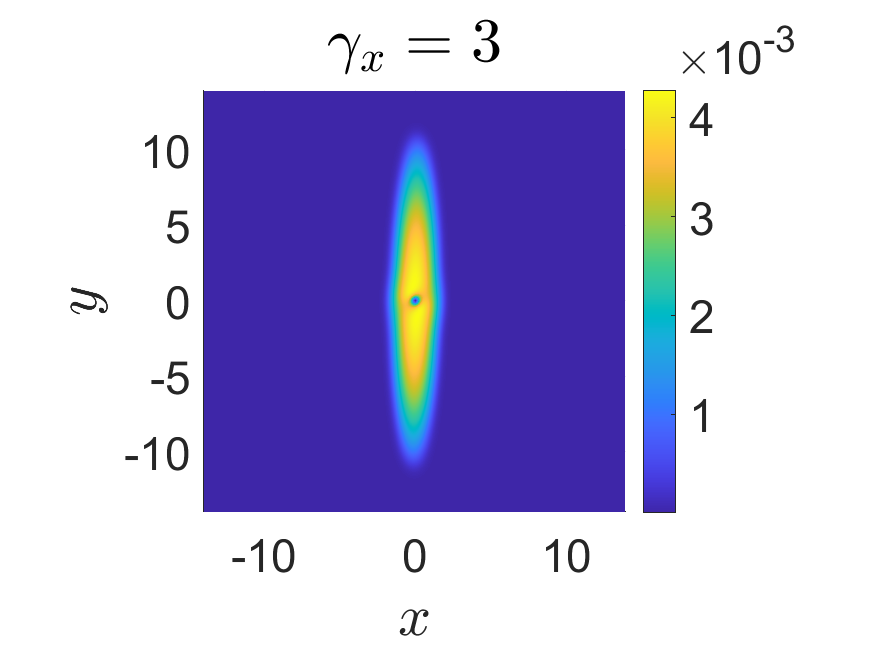,height=2.8cm,width=3.7cm,angle=0}}

\caption{
Contour plots of $|\phi_g(x, y, z=0)|^2$ for \textbf{Case I} (top two rows) and  \textbf{Case II} (bottom two rows) in \textbf{Example \ref{ex_withDDI}}.}\label{vortice_DDI}
\end{figure}

Fig. \ref{vortice_DDI} shows the contour plots of $|\phi_g(x, y, z=0)|^2$ for all cases respectively.
Clearly, the dipole orientation significantly influences the ground state patterns.
The arrangement of the vortices mainly depends on $\bm{n}_\perp$, the $x$-$y$ plane component of the dipole orientation $\bm{n}$.
When $\lambda<0$, the vortices tend to align along the direction that is perpendicular to $\bm{n}_\perp$ (\textbf{Case I});
when $\lambda>0$, the vortices tend to align along with $\bm{n}_\perp$ (\textbf{Case II}).
As for $\bm{n}=(0, 0, 1)^T$,  the distribution of vortices is similar to that in the case of non-dipolar BEC (cf. Fig. \ref{vortice_noDDI}), but more vortices are induced with $\lambda < 0$.

\section{Conclusion}\label{sec:Conclu}
We proposed an accurate and efficient numerical method to compute the ground states of 3D rotating dipolar Bose-Einstein condensates under strongly anisotropic traps.
To address challenges posed by the fast rotation and strongly anisotropic traps, we incorporate the anisotropic truncated kernel method (ATKM)
for dipolar potential evaluation into the preconditioned conjugate gradient method (PCG).
The proposed method is spectrally accurate and efficient, and its memory requirement is  independent of the anisotropy strength.
We presented extensive numerical results to confirm the accuracy and efficiency, together with applications
to the impacts of different model parameters on critical rotational frequency, energies and chemical potential.
In addition, we observed bent vortices (U-shape) in dipolar BEC and carried out a careful and comprehensive numerical computation.
Novel ground state patterns were identified numerically, and the influences of anisotropy strength and dipolar potential on the number and arrangement of vortices
were investigated. Our method can be extended to more general systems, such as multi-component and droplet dipolar BEC,
with appropriate modifications and we shall report them in future work.

\end{document}